\documentstyle[12pt]{article}

\newtheorem{exam}{{\bf Example}}[section]
\newtheorem{lemma}{{\bf Lemma}}[section]
\newtheorem{theo}{{\bf Theorem}}
\newtheorem{prop}{{\bf Proposition}}[section]

\newtheorem{defn}{{\bf Definition}}[section]
\newtheorem{remark}{{\bf Remark}}[section]
\newtheorem{conj}{{\bf Conjecture}}

\font\bbb=msbm10 scaled\magstep1

\newcommand{\RR}{\mbox{\bbb R}}

\newcommand{\XB}{X \hspace{-2.7mm}^{^{\mbox{\bf --}}}}

\newcommand{\La}{\Lambda}

\begin{document}

\title{\Large \bf Lower bound theorem for normal pseudomanifolds}
\author{{\bf Bhaskar Bagchi}$^{\rm a}$, {\bf Basudeb Datta}$^{\rm
\,b,}$\footnote{Corresponding author. \vspace{2mm} \newline \vspace{2mm}
\mbox{} \hspace{3mm} {\em E-mail addresses:} bbagchi@isibang.ac.in
(B. Bagchi), dattab@math.iisc.ernet.in (B. Datta).} $^{, 1}$
}

\date{}

\maketitle

\vspace{-4mm}

\noindent {\footnotesize  $^{\rm a}$Theoretical Statistics and
Mathematics Unit, Indian Statistical Institute,  Bangalore
560\,059, India

\smallskip

\noindent $^{\rm b}$Department of Mathematics, Indian Institute
of Science, Bangalore 560\,012,  India}

\footnotetext[1]{Partially supported by DST (Grant:
SR/S4/MS-272/05) and by UGC-SAP/DSA-IV.
 }

\begin{center}

\date{Appeared in `{\bf Expositiones Mathematicae}'}

\end{center}

\hrule

\medskip

\noindent {\bf Abstract}

{\small In this paper we present a
self-contained combinatorial proof of the lower bound theorem for
normal pseudomanifolds, including a treatment of the cases of
equality in this theorem. We also discuss McMullen and Walkup's
generalised lower bound conjecture for triangulated spheres in
the context of the lower bound theorem. Finally, we pose a new
lower bound conjecture for non-simply connected triangulated
manifolds. }

\bigskip

{\small


\medskip


}

\bigskip

\hrule

\section{Introduction}

The lower bound theorem (LBT) provides the best possible lower
bound for the number of faces of each dimension (in terms of the
dimension and the number of vertices) for any normal
pseudomanifold. When the dimension is at least three, equality
holds precisely for stacked spheres. (This is Theorem \ref{LBT}
in Section 8 below.)

Walkup, Barnette, Klee, Gromov, Kalai and Tay proved various
special cases of the LBT, with Tay providing the first proof in
the entire class of normal pseudomanifolds (cf. \cite{ba3, gv,
ka, kl, ta, wa}). However, Tay's proof rests on Kalai's, and that
in turn depends on the theory of rigidity of frameworks.

Kalai showed in \cite{ka} that for $d \geq 3$, the edge graph of
any connected triangulated $d$-manifold without boundary is
``generically $(d + 1)$-rigid'' in the sense of rigidity of
frameworks. Namely, a particular embedding of a graph in the
$(d+1)$-dimensional Euclidean space is {\em rigid} if it can't be
moved to a nearby embedding without distorting the edge-lengths
(except trivially by bodily moving the entire embedded graph by
applying a rigid motion of the ambient space). A graph is {\em
generically $(d + 1)$-rigid} if the set of its rigid embeddings
in $(d+1)$-space is a dense open subspace in the space of all its
embeddings. The LBT for triangulated manifolds without boundary
is an immediate consequence of Kalai's rigidity theorem. Kalai
also used these ideas to settle the equality case of LBT.
Actually, he proved this theorem in the somewhat larger class of
normal pseudomanifolds whose two-dimensional links are spheres.
In \cite{ta}, Tay showed that Kalai's argument extends almost
effortlessly to the class of all normal pseudomanifolds. This
class has the advantage of being closed under taking links, so
that an induction on dimension is facilitated. Further, the so
called M-P-W reduction (after McMullen, Perles and Walkup) works
in a link-closed class of pseudomanifolds and this reduces the
proof of the general LBT to proving the lower bound only for the
number of edges.

The interesting application of the LBT found in \cite{bd8} led us
to take a close look at Kalai's proof. However, we found it
difficult to follow Kalai's proof in its totality because of our
lack of familiarity with the rigidity theory of frameworks, which
in turn is heavily dependent on analytic considerations that seem
foreign to the questions at hand. We have reasons to suspect that
many experts in Combinatorial Topology share our desire to see a
self-contained combinatorial proof of this fundamental result of
Kalai. For instance, in a relatively recent paper \cite{bb}, Blind
and Blind present a combinatorial proof of the LBT in the class of
polytopal spheres, even though much more general versions were
available. These authors motivate their paper by stating that
``no elementary proof of the LBT including the case of equality
is known so far". One objective of this paper is to rectify this
situation. It may be noted that Blind and Blind use the notion of
shelling to prove the LBT for polytopal spheres. Shelling orders
do not exist in general triangulated spheres (let alone normal
pseudomanifolds), so that the proof presented here is of
necessity very different.

A pointer to a combinatorial proof of LBT for triangulated closed
manifolds was given by Gromov in \cite[pages 211--212]{gv}. There
he introduced a combinatorial analogue of rigidity (which we call
{\em Gromov-rigidity}, or simply {\em rigidity} in this paper) and
sketched an induction argument on the dimension to show that
triangulated $d$-manifolds without boundary are $(d+1)$-rigid in
his sense for $d\geq 2$. However, there was an error at the
starting point $d=2$ of his argument. Reportedly, Connelly and
Whiteley filled this gap, but it seems that their work remained
unpublished. In \cite{ta}, Tay gave a proof of Gromov 3-rigidity
of 2-manifolds. Here we present an independent proof of this
result, based on the notion of generalised bistellar moves
introduced below. It is easy to see that if all the vertex-links
of a $d$-pseudomanifold are Gromov $d$-rigid, then the
$d$-pseudomanifold is $(d+1)$-rigid in the sense of Gromov.
Therefore, $(d+1)$-rigidity of $d$-dimensional normal
pseudomanifolds follows. Now, it is an easy consequence of
Gromov's definition that any $n$-vertex $(d+1)$-rigid simplicial
complex of dimension $d$ satisfies the lower bound $(d+1)n - {d+2
\choose 2}$ on its number of edges, as predicted by LBT. However,
Gromov himself never considered the case of equality in LBT. Here
we refine Gromov's theory to tackle the case of equality. It may
be pointed out that in the concluding remark of \cite{ka}, Kalai
suggested that it should be possible to prove his theorem using
Gromov's ideas. However, the details of such an elementary
argument were never worked out in the intervening twenty years. It
is true that Tay uses Gromov's definition of rigidity in his
proofs. But, to tackle the case of equality, Tay shows that when
equality holds in LBT for a normal pseudomanifold, it must
actually be a triangulated manifold, so that Kalai's initial
argument (based on rigidity of frameworks) applies.

We should note that the notion of generic rigidity pertains
primarily to graphs and Kalai calls a simplicial complex
generically $q$-rigid if its edge graph is generically $q$-rigid.
On the other hand, Gromov's definition pertains to simplicial
complexes. For this reason, it is not possible to compare these
two notions in general. However, such a comparison is possible
when the dimension $d$ of the simplicial complex is $\geq q - 1$
(and we are interested in the case $d = q - 1$). In these cases,
Gromov's notion of rigidity is weaker than the notion of generic
rigidity. From the theory of rigidity of frameworks, it is known
that if an $n$-vertex graph $G$ is minimally generically
$q$-rigid (i.e., $G$ is generically $q$-rigid but no proper
spanning subgraph of $G$ is generically $q$-rigid) then either
$G$ is a complete graph on at most $q + 1$ vertices, or else $G$
has $n \geq q + 1$ vertices and has exactly $nq - {q + 1 \choose
2}$ edges, and any induced subgraph of $G$ (say, with $p \geq q$
vertices) has at most $pq - {q + 1 \choose 2}$ edges (cf.
\cite{gss}. By a theorem of Laman, this fact characterizes
minimally generically $q$-rigid graphs for $q\leq 2$). Using this
result, it is easy to deduce that generic $q$-rigidity (of the
edge graph) implies Gromov's $q$-rigidity for any simplicial
complex of dimension $\geq q-1$.

Apart from the pedagogic/esthetic reason for providing an
elementary proof of the LBT for normal pseudomanifolds (surely an
elementary statement deserves an elementary proof!), we also hope
that the arguments developed here should extend to yield a proof
of the generalised lower bound conjecture (GLBC) for triangulated
spheres. Stanley \cite{st1} proved this conjecture for polytopal
spheres using heavy algebraic tools, but the general case of this
conjecture due to McMullen and Walkup \cite{mw} remains unproved.
Even in Stanley's result, the characterisation of the equality
case remains to be done.

This paper is organized as follows. In Section 2, we give the
preliminary definitions, including an explanation of most of the
technical terms used in this introduction. In the next four
sections, we develop the necessary tools for our proofs. Section
3 provides a combinatorial version of the topological operations
of cutting or pasting handles and of connected sums. These
combinatorial operations were introduced by Walkup in \cite{wa}.
However, the precise combinatorics of these operations was never
worked out. Section 4 introduces the main actors in the game of
LBT's, namely stacked spheres and stacked balls. We also present
some elementary but useful results on these objects. These are
mostly well known, at least to experts. In Section 5, we
introduce the notion of generalized bistellar moves (GBM) and
establish their elementary properties. As the name suggests, this
is a generalization of the usual notion of bistellar moves. It is
also shown that any $n$-vertex triangulated 2-sphere (with $n >
4$) is obtained from an $(n - 1)$-vertex triangulated 2-sphere by
a GBM. More generally, we show that any triangulated orientable
2-manifold (without boundary) $X$ is either the connected sum of
two smaller objects of the same sort, or it is obtained from a
similar object of smaller genus by pasting a handle, or else it
may be obtained by a GBM from a triangulation $\widetilde{X}$ of
the same manifold using one less vertex. (We wonder if similar
results are true for triangulated 3-manifolds.) These results for
triangulated 2-manifolds without boundary are used to give an
inductive proof of their Gromov 3-rigidity in Section 7. Section
6 contains the general theory of Gromov-rigidity, including a
careful treatment of the minimal situations. In Section 7, we
prove the Gromov $(d + 1)$-rigidity of normal $d$-pseudomanifolds,
and show that for $d > 2$ the minimally Gromov $(d + 1)$-rigid
normal pseudomanifolds are precisely the stacked $d$-spheres.
This is Theorem \ref{GRT}, the main result of this paper. As
already indicated, the proof is an induction on $d$. Cutting
handles plays an important role here. In Section 8, we describe
the M-P-W reduction and use it to present the routine deduction of
the LBT for normal pseudomanifolds from Theorem \ref{GRT}. In the
concluding section, we state and discuss the GLBC in a form which
brings out its similarity with the LBT (which is the case $k=1$
of the GLBC). Included in this section is a discussion of the
$k$-stacked spheres which are expected to play a role in the GLBC
similar to the role played by the stacked spheres in LBT. We
conclude by posing a new lower bound conjecture for non-simply
connected triangulated manifolds.


\section{Preliminaries}

Recall that a {\em simplicial complex} is a set of finite sets
such that every subset of an element is also an element. For $i
\geq 0$, an element of size $i + 1$ is called a {\em face} of
dimension $i$ (or an {\em $i$-face)} of the complex. By
convention, the empty set is a face of dimension $-1$. All
simplicial complexes which appear in this paper are finite. The
dimension of a simplicial complex $X$ (denoted by $\dim(X)$) is
by definition the maximum of the dimensions of its faces. The
1-dimensional faces of a simplicial complex are also called the
{\em edges} of the complex. $V(X)$ denotes the set of vertices of
a complex $X$ and is called the {\em vertex-set} of $X$.

For a simplicial complex $X$, $|X|$ is the set of all functions
$f: V(X) \to [0, 1]$ such that $\sum_{v \in V (X)} f(v) = 1$ and
${\rm support}(f) := \{ v \in V (X) : f (v) \neq 0 \}$ is a face
of $X$. (Such a function $f$ may be thought of as a convex
combination of the Dirac delta-functions $\delta_x$ as $x$ ranges
over the face ${\rm support}(f)$.) As a subset of the topological
space $[0, 1]^{V(X)}$, $|X|$ inherits the subspace topology. The
topological space $|X|$ thus obtained is called the {\em
geometric carrier} of $X$. If $|X|$ is a manifold (with or
without boundary) then $X$ is said to be a {\em triangulated
manifold}, or a {\em triangulation} of the manifold $|X|$.

A {\em graph} is a simplicial complex of dimension at most $1$. A
set of vertices of a graph $G$ is said to be a {\em clique} of
$G$ if any two of these vertices are adjacent in $G$ (i.e., form
an edge of $G$). For a general simplicial complex $X$, the {\em
edge graph} (or {\em $1$-skeleton}) $G(X)$ of $X$ is the
subcomplex of $X$ consisting of all its faces of dimensions $\leq
 1$. (More generally, for $0\leq k\leq \dim(X)$, the {\em
$k$-skeleton} ${\rm skel}_k(X)$ of $X$ is the subcomplex
consisting of all the faces of $X$ of dimension $\leq k$.) Notice
that each face of $X$ is a clique in the graph $G(X)$.

If $X$, $Y$ are two simplicial complexes with disjoint vertex
sets, then their {\em join} $X \ast Y$ is the simplicial complex
whose faces are the (disjoint) unions of faces of $X$ with faces
of $Y$. In particular, if $X$ consists of a single vertex $x$,
then we write $x \ast Y$ for $X\ast Y$. The complex $x\ast Y$ is
called the {\em cone} over $Y$ (with cone-vertex $x$).

If $Y$ is a subcomplex of a simplicial complex $X$ and $Y$
consists of all the faces of $X$ contained in $V(Y)$, then we say
that $Y$ is an {\em induced subcomplex} of $X$. If $A \subseteq
V(X)$, then the induced subcomplex of $X$ with the vertex-set $A$
is denoted by $X[A]$. If $\alpha$ is a $k$-face of $X$, then the
{\em closure} $\overline{\alpha}$ of $\alpha$ is the induced
subcomplex $X[\alpha]$. Notice that $\overline{\alpha}$ consists
of all the subsets of $\alpha$. Thus $\overline{\alpha}$ is a
triangulation of the $k$-ball and is also denoted by
$B^{\,k}_{k+1}(\alpha)$.

If $V(X) = A \sqcup B$ is the disjoint union of two subsets $A$
and $B$, then the induced subcomplexes $X[A]$ and $X[B]$ are said
to be {\em simplicial complements} of each other. If $Y$ is an
induced subcomplex of $X$, then the simplicial complement of $Y$
is denoted by $C(Y, X)$. For a face $\alpha$ of $X$, the
simplicial complement $C(\overline{\alpha}, X)$ is called the
{\em antistar} of $\alpha$, and is denoted by ${\rm
ast}(\alpha)$. Thus, ${\rm ast}(\alpha)$ is the subcomplex of $X$
consisting of all faces disjoint from $\alpha$. The {\em link} of
$\alpha$ in $X$, denoted by ${\rm lk}(\alpha)$ (or ${\rm
lk}_X(\alpha)$) is the subcomplex of ${\rm ast}(\alpha)$
consisting of all faces $\beta$ such that $\alpha \sqcup \beta \in
X$. For a vertex $v$ of $X$, the cone $v\ast {\rm lk}_X(v)$ is
called the {\em star} of $v$ in $X$ and is denoted by ${\rm
star}(v)$ (or ${\rm star}_X(v)$).

A $d$-dimensional simplicial complex $X$ is said to be {\em pure}
if all the maximal faces of $X$ have dimension $d$. The maximal
faces in a pure simplicial complex are called its {\em facets}.
The {\em facet graph} $\La (X)$ of a pure $d$-dimensional
simplicial complex $X$ is the graph whose vertices are the facets
of $X$, two such vertices being adjacent in $\La (X)$ if the
corresponding facets intersect in a $(d-1)$-face.

A simplicial complex $X$ is said to be {\em connected} if $|X|$
is connected. Notice that $X$ is connected if and only if its
edge graph $G(X)$ is connected (i.e., any two vertices of $X$ are
the end vertices of a path in $G(X)$). A pure simplicial complex
$X$ is said to be {\em strongly connected} if its facet graph
$\La(X)$ is connected. The connected components of $X$ are the
maximal connected subcomplex of $X$. The strong components of $X$
are the maximal pure subcomplexes of dimension $d = \dim(X)$ which
are strongly connected. Notice that the connected components are
vertex-disjoint, while the strong components may have faces of
codimension two or more in common.

For $d \geq 1$, a $d$-dimensional pure simplicial complex is said
to be a {\em weak pseudomanifold with boundary} if each $(d -
1)$-face is in at most two facets, and it has a $(d-1)$-face
contained in only one facet. A $d$-dimensional pure simplicial
complex is said to be a {\em weak pseudomanifold without
boundary} (or simply {\em weak pseudomanifold}) if each $(d -
1)$-face is in exactly two facets. If $X$ is a $d$-dimensional
weak pseudomanifold with boundary then its boundary $\partial X$
is defined to be the $(d - 1)$-dimensional pure simplicial
complex whose facets are those $(d-1)$-faces of $X$ which are in
unique facets of $X$. Clearly, the link of a face in a weak
pseudomanifold is a weak pseudomanifold.

A {\em pseudomanifold} (respectively {\em pseudomanifold with
boundary}) is a strongly connected weak pseudomanifold
(respectively weak pseudomanifold with  boundary). A
$d$-dimensional weak pseudomanifold (respectively weak
pseudomanifold with boundary) is called a {\em normal
pseudomanifold} (respectively {\em normal pseudomanifold with
boundary}) if each face of dimension $\leq d-2$ has a connected
link. Since we include the empty set as a face, a normal
pseudomanifold is necessarily connected. But we actually have\,:

\begin{lemma}$\!\!\!${\bf .} \label{LPRE1}
Every normal pseudomanifold $($respectively, normal pseudomanifold
with boundary$)$ is a pseudomanifold $($respectively
pseudomanifold with boundary$)$.
\end{lemma}

\noindent {\bf Proof.} Let $X$ be a normal pseudomanifold of
dimension $d\geq 1$. We have to show that its facet graph
$\La(X)$ is connected. If not, choose two facets $\sigma_1$,
$\sigma_2$ from different components of $\La(X)$ for which
$\dim(\sigma_1 \cap \sigma_2)$ is maximum. Then $\dim(\sigma_1
\cap \sigma_2) \leq d-2$ but ${\rm lk}(\sigma_1 \cap \sigma_2)$
is disconnected, a contradiction.  \hfill $\Box$

\bigskip

From the definitions, it is clear that any $d$-dimensional weak
pseudomanifold (respectively weak pseudomanifold with boundary)
has at least $d+2$ (respectively $d + 1$) vertices, with equality
if and only if it is the simplicial complex $S^{\,d}_{d+2}$
(respectively $B^{\,d}_{d+1}$) whose faces are all the proper
subsets of a set of size $d+2$ (respectively, all subsets of a
set of size $d+1$). Clearly, $S^{\,d}_{d+2}$ and $B^{\,d}_{d+1}$
triangulate the $d$-sphere and the $d$-ball, respectively. They
are called the {\em standard $d$-sphere} and the {\em standard
$d$-ball} respectively.

A simplicial complex $X$ is called a {\em combinatorial
$d$-sphere} (respectively, {\em combinatorial $d$-ball\,}) if
$|X|$ (with the induced pl structure from $X$) is pl homeomorphic
to $|S^{\hspace{.2mm}d}_{d + 2}|$ (respectively,
$|B^{\hspace{.1mm}d}_{d + 1}|$).

If $\alpha$ is a face of a simplicial complex $X$, then the number
of vertices in ${\rm lk}_X(\alpha)$ is called the {\em degree} of
$\alpha$ in $X$ and is denoted by $\deg_X(\alpha)$ (or
$\deg(\alpha)$). So, the degree of a vertex $v$ in $X$ is the same
as the degree of $v$ in the edge graph $G(X)$. Since the link of
an $i$-face $\alpha$ in a $d$-dimensional weak pseudomanifold $X$
without boundary is a $(d - i - 1)$-dimensional weak
pseudomanifold, it follows that $\deg_X(\alpha) \geq d - i + 1$,
with equality only if ${\rm lk}_X(\alpha)$ is the standard sphere
$S^{\,d - i - 1}_{d - i + 1}$.

If $X$ is a $d$-dimensional simplicial complex then, for $0\leq j
\leq d$, the number of its $j$-faces is denoted by $f_j =
f_j(X)$. The vector $(f_0, \dots, f_d)$ is called the {\em
face-vector} of $X$ and the number $\chi(X) := \sum_{i=0}^{d}
(-1)^i f_i$ is called the {\em Euler characteristic} of $X$. As
is well known, $\chi(X)$ is a topological invariant, i.e., it
depends only on the homeomorphic type of $|X|$.

\section{Cutting and pasting handles}

\begin{defn}$\!\!\!${\bf .} \label{DAB}
{\rm Let $\sigma_1$, $\sigma_2$ be two facets in a pure
simplicial complex $X$. Let $\psi : \sigma_1 \to \sigma_2$ be a
bijection. We shall say that $\psi$ is {\em admissible} if
($\psi$ is a bijection and) the distance between $x$ and
$\psi(x)$ in the edge graph of $X$ is $\geq 3$ for each $x\in
\sigma_1$ (i.e., if every path in the edge graph joining $x$ to
$\psi(x)$ has length $\geq 3$). Notice that if $\sigma_1$,
$\sigma_2$  are from different connected components of $X$ then
any bijection between them is admissible. Also note that, in
general, for the existence of an admissible map $\psi : \sigma_1
\to \sigma_2$, the facets $\sigma_1$ and $\sigma_2$ must be
disjoint.}
\end{defn}

\begin{defn}$\!\!\!${\bf .} \label{DEHA}
{\rm Let $X$ be a weak pseudomanifold with disjoint facets
$\sigma_1$, $\sigma_2$. Let $\psi \colon \sigma_1\to \sigma_2$ be
an admissible bijection. Let $X^{\psi}$ denote the weak
pseudomanifold obtained from $X \setminus \{\sigma_1, \sigma_2\}$
by identifying $x$ with $\psi(x)$ for each $x\in \sigma_1$. Then
$X^{\psi}$ is said to be obtained from $X$ by an {\em elementary
handle addition}. If $X_1$, $X_2$ are two $d$-dimensional weak
pseudomanifolds with disjoint vertex-sets, $\sigma_i$ a facet of
$X_i$ ($i=1, 2$) and $\psi \colon \sigma_1 \to \sigma_2$ any
bijection, then $(X_1\sqcup X_2)^{\psi}$ is called an {\em
elementary connected sum} of $X_1$  and $X_2$, and is denoted by
$X_1 \#_{\psi} X_2$ (or simply by $X_1\# X_2$). Note that the
combinatorial type of $X_1 \#_{\psi} X_2$ depends on the choice
of the bijection $\psi$. However, when $X_1$, $X_2$ are connected
triangulated $d$-manifolds, $|X_1 \#_{\psi} X_2|$ is the
topological connected sum of $|X_1|$ and $|X_2|$ (taken with
appropriate orientations). Thus, $X_1 \#_{\psi} X_2$ is a
triangulated manifold whenever $X_1$, $X_2$ are triangulated
$d$-manifolds. }
\end{defn}

\begin{lemma}$\!\!\!${\bf .} \label{sd-s}
Let $N$ be a $(d- 1)$-dimensional induced subcomplex of a
$d$-dimensional simplicial complex $M$. If both $M$ and $N$ are
normal pseudomanifolds then \begin{enumerate}
\item[$(a)$] for any vertex $u$ of $N$ and any vertex $v$ of the
simplicial complement $C(N, M)$, there is a path $P$ $($in $M)$
joining $u$ to $v$ such that $u$ is the only vertex in $P\cap N$,
and
\item[$(b)$] the simplicial complement $C(N, M)$ has at most two
connected components.
\end{enumerate}
\end{lemma}

\noindent {\bf Proof.}  Part $(a)$ is trivial if $d=1$ (in which
case, $N = S^{\hspace{.2mm}0}_2$ and $M = S^{\hspace{.2mm}1}_n$).
So, assume $d > 1$ and we have the result for smaller dimensions.
Clearly, there is a path $P$ (in the edge graph of $M$) joining
$u$ to $v$ such that $P = x_1x_2 \cdots x_ky_1 \cdots y_l$ where
$x_1 = u$, $y_l = v$ and $x_i$'s are the only vertices of $P$
from $N$. Choose $k$ to be the smallest possible. We claim that $k
=1$, so that the result follows. If not, then $x_{k-1} \in {\rm
lk}_N(x_k) \subset {\rm lk}_M(x_k)$ and $y_1 \in C({\rm
lk}_N(x_k), {\rm lk}_M(x_k))$. Then, by induction hypothesis,
there is a path $Q$ in ${\rm lk}_M(x_k)$ joining $x_{k-1}$ and
$y_1$ in which $x_{k-1}$ is the only vertex from ${\rm
lk}_N(x_k)$. Replacing the part $x_{k-1}x_ky_1$ of $P$ by the
path $Q$, we get a path $P^{\hspace{.2mm} \prime}$ from $u$ to
$v$ where only the first $k-1$ vertices of $P^{\hspace{.2mm}
\prime}$ are from $N$. This contradicts the choice of $k$.

The proof of Part $(b)$ is also by induction on the dimension $d$.
The result is trivial for $d=1$. For $d > 1$, fix a vertex $u$ of
$N$. By induction hypothesis, $C({\rm lk}_N(u), {\rm lk}_M(u))$
has at most two connected components. By Part $(a)$ of this
lemma, every vertex $v$ of $C(N, M)$ is joined by a path in $C(N,
M)$ to a vertex in one of these components. Hence the result.
\hfill $\Box$

\bigskip

Let $N$ be an induced subcomplex of a simplicial complex $M$. One
says that $N$ is {\em two-sided} in $M$ if $|N|$ has a (tubular)
neighbourhood in $|M|$ homeomorphic to $|N| \times [-1, 1]$ such
that the image of $|N|$ (under this homeomorphism) is $|N|\times
\{0\}$.

\begin{lemma}$\!\!\!${\bf .} \label{m-sts}
Let $M$ be a normal pseudomanifold of dimension $d \geq 2$ and
$A$ be a set of vertices of $M$ such that the induced subcomplex
$M[A]$ of $M$ on $A$ is a $(d - 1)$-dimensional normal
pseudomanifold. Let $G$ be the graph whose vertices are the edges
of $M$ with exactly one end in $A$, two such vertices being
adjacent in $G$ if the union of the corresponding edges is a
$2$-face of $M$. Then $G$ has at most two connected components.
If, further, $M[A]$ is two-sided in $M$ then $G$ has exactly two
connected components.
\end{lemma}

\noindent {\bf Proof.} Let $E = V(G)$ be the set of edges of $M$
with exactly one end in $A$. For $x\in A$, set $E_x = \{e\in E :
x\in e\}$, and let $G_x = G[E_x]$ be the induced subgraph of $G$
on $E_x$. Note that $G_x$ is isomorphic to the edge graph  of
$C({\rm lk}_{M[A]}(x), {\rm lk}_M(x))$. Therefore, by Lemma
\ref{sd-s} $(b)$, $G_x$ has at most two components for each $x\in
A$. Also, for an edge $xy$ in $M[A]$, there is a $d$-face
$\sigma$ of $M$ such that $xy$ is in $\sigma$. Since the induced
complex $M[A]$ is $(d - 1)$-dimensional, there is a vertex $u \in
\sigma \setminus A$. Then $e_1 = xu \in E_x$ and $e_2 = yu \in
E_y$ are adjacent in $G$. Thus, if $x$, $y$ are adjacent vertices
in $M[A]$ then there is an edge of $G$ between $E_x$ and $E_y$.
Since $M[A]$ is connected and $V(G) = \cup_{x \in A} E_x$, it
follows that $G$ has at most two connected components.

Now suppose $S = M[A]$ is two-sided in $M$. Let $U$ be a tubular
neighbourhood of $|S|$ in $|M|$ such that $U \setminus |S|$ has
two components, say $U^{+}$ and $U^{-}$. Since $|S|$ is compact,
we can choose $U$ sufficiently small so that $U$ does not contain
any vertex from $V(M) \setminus A$. Then, for $e \in E$, $|e|$
meets either $U^{+}$ or $U^{-}$ but not both. Put $E^{\pm} = \{e
\in E : |e| \cap U^{\pm} \neq \emptyset\}$. Then no element of
$E^{+}$ is adjacent in $G$ with any element of $E^{-}$. From the
previous argument, one sees that each $x\in A$ is in an edge from
$E^{+}$ and in an edge from $E^{-}$. Thus, both $E^{+}$ and
$E^{-}$ are non-empty. So, $G$ is disconnected. \hfill $\Box$

\begin{lemma}$\!\!\!${\bf .} \label{LEHD}
Let $X$ be a normal $d$-pseudomanifold with an induced two-sided
standard $(d - 1)$-sphere $S$. Then there is a $d$-dimensional
weak pseudomanifold $\widetilde{X}$ such that $X$ is obtained
from $\widetilde{X}$ by elementary handle addition. Further,
\vspace{-1mm}
\begin{enumerate}
\item[$(a)$]the connected components of $\widetilde{X}$ are normal
$d$-pseudomanifolds, \vspace{-1mm}
\item[$(b)$] $\widetilde{X}$ has at most two connected components,
\vspace{-1mm}
\item[$(c)$] if $\widetilde{X}$ is not connected, then $X = Y_1 \#
Y_2$, where $Y_1$, $Y_2$ are the connected components of
$\widetilde{X}$, and \vspace{-1mm}
\item[$(d)$] if \, $C(S, X)$ is connected then $\widetilde{X}$ is
connected.
\end{enumerate}
\end{lemma}

\noindent {\bf Proof.} As above, let $E$ be the set of all edges
of $X$ with exactly one end in $S$. Let $E^{+}$ and $E^{-}$ be
the connected components of the graph $G$ (with vertex-set $E$)
defined above (cf. Lemma \ref{m-sts}). Notice that if a facet
$\sigma$ intersects $V(S)$ then $\sigma$ contains edges from $E$,
and the graph $G$ induces a connected subgraph on the set
$E_{\sigma} = \{e\in E : e\subseteq \sigma\}$. (Indeed, this
subgraph is the line graph of a complete bipartite graph.)
Consequently, either $E_{\sigma} \subseteq E^{+}$ or $E_{\sigma}
\subseteq E^{-}$. Accordingly, we say that the facet $\sigma$ is
positive or negative (relative to $S$). If a facet $\sigma$ of
$X$ does not intersect $V(S)$ then we shall say that $\sigma$ is
a neutral facet.

Let $V(S) = W$ and $V(X) \setminus V(S) = U$. Take two disjoint
sets $W^{+}$ and $W^{-}$, both disjoint from $U$, together with
two bijections $f_{\pm} \colon W \to W^{\pm}$. We define a pure
simplicial complex $\widetilde{X}$ as follows. The vertex-set of
$\widetilde{X}$ is $U \sqcup W^{+} \sqcup W^{-}$. The facets of
$\widetilde{X}$  are: (i) $W^{+}$, $W^{-}$, (ii) all the neutral
facets of $X$, (iii) for each positive facet $\sigma$ of $X$, the
set $\widetilde{\sigma} := (\sigma \cap U) \sqcup f_{+}(\sigma
\cap W)$, and (iv) for each negative facet $\tau$ of $X$, the set
$\widetilde{\tau} := (\tau \cap U) \sqcup f_{\!-}(\tau \cap W)$.
Clearly, $\widetilde{X}$ is a weak pseudomanifold. Let $\psi =
f_{-} \circ f_{+}^{-1} \colon W^{+} \to W^{-}$. It is easy to see
that $\psi$ is admissible and $X = (\widetilde{X})^{\psi}$.

Since the links of faces of dimension up to $d - 2$ in $X$ are
connected, it follows that the links of faces of dimension up to
$d - 2$ in $\widetilde{X}$ are connected. This proves $(a)$.

As $X$ is connected, choosing two vertices $f_{\pm}(x_0) \in
W^{\pm}$ of $\widetilde{X}$, one sees that each vertex of
$\widetilde{X}$ is joined by a path in the edge graph  of
$\widetilde{X}$ to either $f_{+}(x_0)$ or $f_{-}(x_0)$. Hence
$\widetilde{X}$ has at most two components. This proves $(b)$.
This arguments also shows that when $\widetilde{X}$ is
disconnected, $W^{+}$ and $W^{-}$ are facets in different
components of $\widetilde{X}$. Hence $(c)$ follows.

Observe that $C(S, X) = C(W^{+} \sqcup W^{-}, \widetilde{X})$.
Assume that $C(S, X)$ is connected. Now, for any $(d-1)$-simplex
$\tau \subseteq W^{+}$, there is a vertex $x$ in $C(S, X)$ such
that $\tau\cup\{x\}$ is a facet of $\widetilde{X}$. So, $C(S, X)$
and $W^{+}$ are in the same connected component of
$\widetilde{X}$. Similarly, $C(S, X)$ and $W^{-}$ are in the same
connected component of $\widetilde{X}$. This proves $(d)$. \hfill
$\Box$

\begin{defn}$\!\!\!${\bf .} \label{DEHD}
{\rm If $S$ is an induced two-sided $S^{\hspace{.35mm}d - 1}_{d +
1}$ in a normal $d$-pseudomanifold $X$, then the pure simplicial
complex $\widetilde{X}$ constructed above is said to be obtained
from $X$ by an} elementary handle deletion {\rm over $S$.}
\end{defn}

\begin{remark}$\!\!\!${\bf .} \label{REHD}
{\rm In Lemma \ref{LEHD}, if $X$ is a triangulated manifold then
it is easy to see that $\widetilde{X}$ is also a triangulated
manifold.}
\end{remark}

\begin{exam}$\!\!\!${\bf .} \label{RP26}
{\rm It is well known that the real projective plane has a unique
$6$-vertex triangulation, denoted by $\RR P^{\hspace{.2mm}2}_6$.
It is obtained from the boundary complex of the icosahedron by
identifying antipodal vertices. The simplicial complement of any
facet in $\RR P^{\hspace{.2mm}2}_6$ is an $S^{\hspace{.2mm}1}_3$.
But, it is not possible to obtain a triangulated 2-manifold $M$
by deleting the handle over this $S^{\hspace{.2mm}1}_3$. Such a
2-manifold would have face vector $(9, 18, 12)$ and hence Euler
characteristic $\chi = 3$. But, arguing as in the proof of Lemma
\ref{LEHD} $(d)$, one can see that $M$ must be connected - and
any connected closed 2-manifold has Euler characteristic $\leq
2$, a contradiction. Thus the hypothesis ``two-sided" in
Definition \ref{DEHD} is essential. Indeed, in this example, the
graph $G$ of Lemma \ref{m-sts} is connected: it is a $9$-gon.}
\end{exam}


\section{Stacked spheres}

Let $X$ be a pure $d$-dimensional weak pseudomanifold and $\sigma$
be a facet of $X$. Take a symbol $v$ outside $V(X)$, and $Y$ be
the pure simplicial complex with vertex set $V(X) \cup \{v\}$
whose facets are facets of $X$ other than $\sigma$ and the $(d +
1)$-sets $\tau \cup \{v\}$ where $\tau$ runs over the $(d -
1)$-faces in $\sigma$. Clearly, $Y$ is a weak pseudomanifold and
$|X|$ and $|Y|$ are homeomorphic topological spaces. This $Y$ is
said to be the weak pseudomanifold obtained from $X$ by {\em
starring} the new vertex $v$ in the facet $\sigma$. (In the
literature, this is also known as the {\em bistellar $0$-move}.)
Notice that the new vertex $v$ is of (minimal) degree $d+1$ in
$Y$. Conversely, let $Y$ be a $d$-dimensional weak pseudomanifold
with a vertex $v$ of degree $d+1$. Let $\sigma = V({\rm
lk}_Y(v))$. If $\sigma$ is not a face of $Y$ (which is
automatically true if $Y$ is a pseudomanifold other than the
standard $d$-sphere $S^{\,d}_{d+2}$) then consider the pure
simplicial complex $X$ with vertex-set $V(Y)\setminus\{v\}$ whose
facets are the facets of $Y$ not containing $v$ and the
$(d+1)$-set $\sigma$. Clearly, $X$ is a weak pseudomanifold. This
$X$ is said to be obtained from $Y$ by {\em collapsing} the
vertex $v$. (This is also called a {\em bistellar $d$-move} in
the literature.) Obviously, the operations of starring a vertex
in a facet and collapsing a vertex of minimal degree are inverses
of each other.

\begin{defn}$\!\!\!${\bf .} \label{DSS1}
{\rm A simplicial complex $X$ is said to be a {\em stacked
$d$-sphere} if $X$ is obtained from the standard $d$-sphere
$S^{\,d}_{d+2}$ by a finite sequence of bistellar 0-moves.
Clearly, any stacked $d$-sphere is a combinatorial $d$-sphere.}
\end{defn}

\begin{lemma}$\!\!\!${\bf .} \label{LSSSB1}
Let $X$ be a triangulated $d$-sphere and $x$ be a vertex of $X$.
If ${\rm lk}_X(x)$ is a triangulated sphere then ${\rm ast}_X(x)$
is a triangulated $d$-ball. In particular, if $X$ is a
combinatorial $d$-sphere then the antistar of every vertex of $X$
is a triangulated ball.
\end{lemma}

\noindent {\bf Proof.} Note that $|{\rm ast}_X(x)|$ is the
closure of a component of $|X|\setminus |{\rm lk}_X(x)|$. Also,
$|{\rm lk}_X(x)|$ has a neighbourhood in $|X|$ which is
homeomorphic to $|{\rm lk}_X(x)|\times [-1, 1]$ via a
homeomorphism mapping $|{\rm lk}_X(x)|$ onto $|{\rm
lk}_X(x)|\times \{0\}$. Therefore, by the generalized
Sch\"{o}nflies theorem (cf. \cite[Theorem 5]{br}), $|{\rm
ast}_X(x)|$ is a $d$-ball. If $X$ is a combinatorial $d$-sphere,
then each vertex link is a triangulated (indeed combinatorial)
sphere, so that this argument applies to each vertex of $X$.
\hfill $\Box$

\begin{defn}$\!\!\!${\bf .} \label{DSBall}
{\rm A {\em stacked $d$-ball} is by definition the antistar of a
vertex in a stacked $d$-sphere. Thus if $X$ is a stacked
$d$-sphere and $x$ is a vertex of $X$, then the simplicial
complex $Y$, whose faces are the faces of $X$ not containing $x$,
is a stacked $d$-ball. Lemma \ref{LSSSB1} implies that stacked
$d$-balls are indeed triangulated balls. It's not hard to see
that they are actually combinatorial balls.}
\end{defn}

From the above discussion, we see that any stacked $d$-sphere is
a triangulation of the $d$-dimensional sphere. Since an
$n$-vertex stacked $d$-sphere is obtained from $S^{\,d}_{d + 2}$
by $(n - d - 2)$ starring and each starring induces ${d+1 \choose
j}$ new $j$-faces and retains all the old $j$-faces for $1\leq j
< d$ (respectively, kills only one old $j$-face for $j=d$), it
follows that it has $(n-d-2){d+1 \choose j} + {d+2 \choose j+1}$
$j$-faces for $1\leq j < d$, and $(n-d-2)d + (d+2)$ facets. On
simplifying, we get\,:

\begin{lemma}$\!\!\!${\bf .} \label{LSSSB2}
The face-vector of any $d$-dimensional stacked sphere satisfies
$$
f_j = \left\{\begin{array}{ll}
 {d+1 \choose j}f_0 - j{d+2 \choose j+1}, & \mbox{if } ~~ 1\leq j<d \\
 df_0 - (d+2)(d-1), & \mbox{if } ~~ j=d.
 \end{array}
 \right.
$$
\end{lemma}

\begin{lemma}$\!\!\!${\bf .} \label{LSSSB3}
Let $X$ be a normal pseudomanifold of dimension $d \geq 2$.
\vspace{-2mm}
\begin{enumerate}
\item[{$(a)$}] If $X \neq S^{\,d}_{d+2}$ then any two vertices of
degree $d + 1$ in $X$ are non-adjacent. \vspace{-2mm}
\item[{$(b)$}] If $X$ is a stacked $d$-sphere then $X$ has at
least two vertices of degree $d + 1$.
\end{enumerate}
\end{lemma}

\noindent {\bf Proof.} Let $x_1$, $x_2$ be two adjacent vertices
of degree $d + 1$ in $X$. Thus, ${\rm lk}(x_1) = S^{\hspace{.2mm}d
- 1}_{d + 1}$, so that all the vertices in $V = V({\rm st}(x_1))$
are adjacent. It follows that $V \setminus \{x_2\}$ is the set of
neighbours of $x_2$. Hence all the facets through $x_2$ are
contained in the $(d + 2)$-set $V$. Since there must be a facet
containing $x_2$ but not containing $x_1$, such a facet must be
$V \setminus \{x_1\}$. Thus, $X$ induces a standard $d$-sphere on
$V$. Since $X$ is a $d$-dimensional normal pseudomanifold, it
follows that $X = S^{\hspace{.2mm}d}_{d + 2}(V)$. This proves Part
$(a)$.

We prove $(b)$ by induction on the number $n$ of vertices of $X$.
If $n = d + 2$ then $X = S^{\hspace{.2mm}d}_{d + 2}$ and the
result is trivial. So assume $n > d + 2$, and the result holds
for all the smaller values of $n$. Since $X$ is a stacked sphere,
$X$ is obtained from an $(n - 1)$-vertex stacked sphere $Y$ by
starring a new vertex $x$ in a facet $\sigma$ of $Y$. Thus, $x$ is
a vertex of degree $d + 1$ in $X$. If $Y$ is the standard
$d$-sphere then the unique vertex $y$ in $V(Y) \setminus \sigma$
is also of degree $d + 1$ in $X$. Otherwise, by induction
hypothesis, $Y$ has at least two vertices of degree $d + 1$, and
since any two of the vertices in $\sigma$ are adjacent in $Y$ -
Part $(a)$ implies that at least one of these degree $d + 1$
vertices of $Y$ is outside $\sigma$. Say $z \not\in \sigma$ is of
degree $d + 1$ in $Y$. Then $z$ (as well as $x$) is a vertex of
degree $d + 1$ in $X$.  \hfill $\Box$

\begin{lemma}$\!\!\!${\bf .} \label{LSSSB4}
Let $X$, $Y$ be $d$-dimensional normal pseudomanifolds. Suppose
$Y$ is obtained from $X$ by starring a new vertex in a facet of
$X$. Then $Y$ is a stacked sphere if and only if $X$ is a stacked
sphere.
\end{lemma}

\noindent {\bf Proof.} The ``if\hspace{.5mm}" part is immediate
from the definition of stacked spheres. We prove the ``only
if\hspace{.5mm}" part by induction on the number $n \geq d + 3$
of vertices of $Y$. The result is trivial for $n = d + 3$. So,
assume $n > d + 3$. Let $Y$ be obtained from $X$ by starring a
vertex $x$ in a facet $\sigma$ of $X$. Suppose $Y$ is a stacked
sphere. Then $Y$ is obtained from some stacked sphere $Z$ by
starring a vertex $y$ in a facet $\tau$ of $Z$. If $x = y$ then
$Z$ is obtained from $Y$ by collapsing $x$, so that $X = Z$ is a
stacked sphere, hence we are done. On the other hand, if $x \neq
y$, then both $x$ and $y$ are of degree $d + 1$ in $Y$, so that
by Lemma \ref{LSSSB3}, $x$ and $y$ are non-adjacent. Therefore,
$x$ is a vertex of degree $d + 1$ in $Z$. Let $W$ be obtained from
$Z$ by collapsing the vertex $x$. By induction hypothesis, $W$ is
a stacked sphere. But, $X$ is obtained from $W$ by starring the
vertex $y$. Hence by the ``if\," part, $X$ is a stacked sphere.
\hfill $\Box$

\begin{lemma}$\!\!\!${\bf .} \label{LSSSB5}
The link of a vertex in a stacked sphere is a stacked sphere.
\end{lemma}

\noindent {\bf Proof.} Let $X$ be a $d$-dimensional stacked
sphere and $v$ be a vertex of $X$. We prove the result by
induction on the number $n$ of vertices of $X$. The result is
trivial for $n = d + 2$. So, assume $n \geq d + 3$ and the result
is true for all stacked spheres on at most $n - 1$ vertices. Let
$X$ be obtained from an $(n - 1)$-vertex stacked sphere $Y$ by
starring a vertex $x$ in a facet $\sigma$ of $Y$. If $v = x$ then
${\rm lk}_X(v)$ is a standard $(d-1)$-sphere and hence is a
stacked sphere. So, assume that $v \neq x$. Since the number of
vertices in $Y$ is $n-1$, by induction hypothesis, ${\rm
lk}_Y(v)$ is a stacked sphere. Clearly, either ${\rm lk}_X(v) =
{\rm lk}_Y(v)$ or ${\rm lk}_X(v)$ is obtained from ${\rm
lk}_Y(v)$ by starring $x$ in a facet of ${\rm lk}_Y(v)$. In either
case, ${\rm lk}_X(v)$ is a stacked sphere. \hfill $\Box$

\begin{lemma}$\!\!\!${\bf .} \label{LSSSB6}
Let $X$ be a stacked $d$-sphere with edge graph $G$ and $d
> 1$. Let $\XB$ denote the simplicial complex whose faces are all
the cliques of $G$. Then $\XB$ is a stacked $(d+ 1)$-ball and $X
= \partial \XB$.
\end{lemma}

\noindent {\bf Proof.} We prove the result by induction on the
number $n$ of vertices of $X$. If $n = d+2$ then $X =
S^{\,d}_{d+2}$ and $\XB = B^{\,d+1}_{d+2}$, so that the result is
obviously true. So assume that $n> d+2$ and the result is true
for $(n-1)$-vertex stacked $d$-spheres. Let $x$ be a vertex of
degree $d+1$ in $X$, and let $X_0$ be the $(n-1)$-vertex stacked
$d$-sphere obtained from $X$ by collapsing the vertex $x$. Note
that, since $d\geq 2$, the edge graph $G_0$ of $X_0$ is the
induced subgraph on the vertex-set $V(G_0) = V(G) \setminus
\{x\}$, and $G$ may be recovered from $G_0$ by adding the vertex
$x$ and making it adjacent to the vertices in a $(d+1)$-clique
$\sigma$ of $G_0$ (which formed a facet of $X_0$, i.e., a
boundary $d$-face of the stacked $(d+1)$-ball ${\XB}_0$). Thus
the simplicial complex $\XB$ is obtained from the stacked
$(d+1)$-ball ${\XB}_0$ by adding the $(d+1)$-face $\tilde{\sigma}
:= \sigma \cup\{x\}$. Since ${\XB}_0$ is a stacked $(d+1)$-ball,
it is the antistar of a (new) vertex $y$ in a stacked
$(d+1)$-sphere $Y_0$ with vertex set $V(X_0)\sqcup \{y\}$. Since
$\sigma$ is a boundary face of ${\XB}_0$, it follows that
$\hat{\sigma} := \sigma\sqcup \{y\}$ is a facet of $Y_0$. Let $Y$
be the $(n+1)$-vertex stacked $(d+1)$-sphere obtained from $Y_0$
by starring the vertex $x$ in the facet $\hat{\sigma}$. Clearly,
$\XB$ is the antistar in $Y$ of the vertex $y$. Therefore, $\XB$
is a stacked $(d+1)$-ball. Now, ${\rm lk}_Y(y)$ is obtained from
${\rm lk}_{Y_0}(y)$ by starring the vertex $x$ in the $d$-face
$\sigma$. Since ${\rm lk}_Y(y) =
\partial \XB$ and ${\rm lk}_{Y_0}(y) = \partial {\XB}_0 = X_0$, it
follows that $\partial \XB$ is obtained from $X_0$ by starring
the vertex $x$ in the facet $\sigma$. That is, $\partial \XB = X$.
This completes the induction and hence proves the lemma.  \hfill
$\Box$

\begin{lemma}$\!\!\!${\bf .} \label{LSSSB7}
Any stacked sphere is uniquely determined by its edge graph.
\end{lemma}

\noindent {\bf Proof.} Let $G$ be the edge-graph of a stacked
$d$-sphere $X$. If $d = 1$ then $X = G$, and there is nothing to
prove. If $d > 1$, then Lemma \ref{LSSSB6} shows that $G$
determines $\XB$ (by definition) and $\XB$ determines $X$ via the
formula $X = \partial \XB$. \hfill $\Box$

\begin{remark}$\!\!\!${\bf .} \label{1-stacked}
{\rm (a) From the definition and Lemma \ref{LSSSB5}, it follows
that the boundary of any stacked ball is a stacked sphere.
Conversely, from Lemma \ref{LSSSB6}, every stacked $d$-sphere $X$
is the boundary of a stacked $(d + 1)$-ball $\XB$ canonically
constructed from $X$ for $d \geq 2$. Indeed, $\XB$ is the unique
triangulated ball such that ${\rm skel}_{d - 1}(X) = {\rm
skel}_{d - 1}({\XB})$. Thus, any stacked sphere is a 1-stacked
sphere as defined in Section 9. \newline (b) Lemma \ref{LSSSB2}
implies that any stacked $d$-ball with $n$ boundary vertices and
$m$ interior vertices has exactly $n + (m - 1)d$ facets. In
particular, if $X$ is an $n$-vertex stacked $d$-sphere, then the
stacked $(d + 1)$-ball $\XB$ constructed above has $n$ boundary
vertices and no interior vertices, so that $X$ has exactly $n - d
- 1$ cliques of size $d + 2$. Of course, this may be directly
verified by induction on $n$.}
\end{remark}

\begin{lemma}$\!\!\!${\bf .} \label{LSSSB8}
Let $X_1$, $X_2$ be $d$-dimensional normal pseudomanifolds. Then
$(a)$ $X_1 \# X_2$ is a triangulated $2$-sphere if and only if
both $X_1$ and $X_2$ are triangulated $2$-spheres; and $(b)$ $X_1
\# X_2$ is a stacked $d$-sphere if and only if both $X_1$, $X_2$
are stacked $d$-spheres.
\end{lemma}

\noindent {\bf Proof.} Let $d = 2$. Then $X_1$, $X_2$ are
connected triangulated 2-manifolds and hence $X_1 \# X_2$ is a
connected triangulated 2-manifold. For $0\leq i\leq 2$, $1\leq
j\leq 2$, let $f_i(X_j)$ denote the number of $i$-faces in $X_j$.
Then, from the definition, $\chi(X_1 \# X_2) = (f_0(X_1) +
f_0(X_2) - 3) - (f_1(X_1) + f_1(X_2) - 3) + (f_2(X_1) + f_2(X_2)
- 2) = \chi(X_1) + \chi(X_2) - 2$. Part $(a)$ now follows from
the fact that the Euler characteristic of a connected closed
2-manifold $M$ is $\leq 2$ and equality holds if and only if $M$
is a 2-sphere.

We prove Part $(b)$ by induction on the number $n \geq d + 3$ of
vertices in $X_1 \# X_2$. If $n = d + 3$ then both $X_1$, $X_2$
must be standard $d$-spheres (hence stacked spheres) and then
$X_1 \# X_2 = S^{\, 0}_{2} \ast S^{\hspace{.2mm}d - 1}_{d + 1}$
is easily seen to be a stacked sphere. So, assume $n
> d + 3$, so that at least one of $X_1$, $X_2$ is not the
standard $d$-sphere. Without loss of generality, say $X_1$ is not
the standard $d$-sphere. Of course, $X = X_1 \# X_2$ is not a
standard $d$-sphere. Let $X$ be obtained from $X_1 \sqcup X_2
\setminus \{\sigma_1, \sigma_2\}$ by identifying a facet
$\sigma_1$ of $X_1$ with a facet $\sigma_2$ of $X_2$ by some
bijection. Then, $\sigma_1 = \sigma_2$ is a clique in the edge
graph of $X$, though it is not a facet of $X$. Notice that a
vertex $x \in V(X_1) \setminus \sigma_1$ is of degree $d + 1$ in
$X_1$ if and only if it is of degree $d + 1$ in $X$. If either
$X_1$ is a stacked sphere or $X$ is a stacked sphere then, by
Lemma \ref{LSSSB3}, such a vertex $x$ exists. Let
$\widetilde{X}_1$ (respectively, $\widetilde{X}$) be obtained
from $X_1$ (respectively, $X$) by collapsing this vertex $x$.
Notice that $\widetilde{X} = \widetilde{X}_1 \# X_2$. Therefore,
by induction hypothesis and Lemma \ref{LSSSB4}, we have: $X$ is a
stacked sphere $\Longleftrightarrow$ $\widetilde{X}$ is a stacked
sphere $\Longleftrightarrow$ both $\widetilde{X}_1$ and $X_2$ are
stacked spheres $\Longleftrightarrow$ both $X_1$ and $X_2$ are
stacked spheres.  \hfill $\Box$

\begin{defn}$\!\!\!${\bf .} \label{DSS2}
{\rm For $d \geq 2$, ${\cal K}(d)$ will denote the family of all
$d$-dimensional normal pseudomanifolds $X$ such that the link of
each vertex of $X$ is a stacked $(d - 1)$-sphere. Since all
stacked spheres are combinatorial spheres, it follows that the
members of ${\cal K}(d)$ are triangulated $d$-manifolds.  }
\end{defn}

\begin{lemma}{\bf (Walkup \cite{wa}).} \label{LSSSB9} Let $X$ be
a normal $d$-pseudomanifold and $\psi \colon \sigma_1 \to
\sigma_2$ be an admissible bijection, where $\sigma_1, \sigma_2$
are facets of $X$. Then $(a)$ $X^{\psi}$ is a triangulated
$3$-manifold if and only if $X$ is a triangulated $3$-manifold;
and $(b)$ $X^{\psi} \in {\cal K}(d)$ if and only if $X \in {\cal
K}(d)$.
\end{lemma}

\noindent {\bf Proof.} For a vertex $v$ of $X$, let $\bar{v}$
denote the corresponding vertex of $X^{\psi}$. Observe that ${\rm
lk}_{X^{\psi}}(\bar{v})$ is isomorphic to ${\rm lk}_{X}(v)$ if $v
\in V(X) \setminus (\sigma_1 \cup \sigma_2)$ and ${\rm
lk}_{X^{\psi}}(\bar{v}) = {\rm lk}_{X}(v) \# {\rm
lk}_{X}(\psi(v))$ if $v \in \sigma_1$. The results now follow
from Lemma \ref{LSSSB8}. \hfill $\Box$

\bigskip

Notice that, Lemma \ref{LSSSB5} says that all stacked $d$-spheres
belong to the class ${\cal K}(d)$. Indeed, we have the following
characterization of stacked spheres of dimension $\geq 4$. This
is essentially a result from Kalai \cite{ka}.

\begin{lemma} \label{LSSSB10} For $d \geq 4$, every member of
${\cal K}(d)$, excepting $S^{\,d}_{d + 2}$, has an $S^{\,d - 1}_{d
+ 1}$ as an induced subcomplex.
\end{lemma}

\noindent {\bf Proof.} Let $X \in {\cal K}(d)$, $X\neq S^{\,d}_{d
+ 2}$. Then $X$ has a vertex of degree $\geq d+2$. Fix such a
vertex $x$, and let $\sigma$ be an interior $(d-1)$-face in the
stacked $d$-ball $\overline{{\rm lk}_X(x})$. (If there was no
such $(d-1)$-face, then we would have $\overline{{\rm lk}_X(x}) =
B^{\,d}_{d+1}$, and hence $\deg(x) = d+1$, contrary to the choice
of $x$.) We claim that $X$ induces an $S^{\,d - 1}_{d + 1}$ on
$\sigma\cup\{x\}$. In other words, the claim is that $\sigma\in
X$.

Choose any vertex $y\in \sigma$, and let $\sigma^{\prime} =
(\sigma \cup \{x\})\setminus\{y\}$. Since ${\rm lk}_X(x)$ and
$\overline{{\rm lk}_X(x})$ have the same $(d - 2)$-skeleton and
$\sigma$ is a $(d - 1)$-face of the latter, it follows that every
proper subset of $\sigma^{\prime} \cup \{y\} = \sigma \cup \{x\}$
which contains $x$ is a face of $X$. Since $d\geq 4$, it follows
in particular that $\sigma^{\prime}$ is a clique of the edge
graph of ${\rm lk}_X(y)$. Hence $\sigma^{\prime} \in
\overline{{\rm lk}_X(y})$. Thus every proper subset of
$\sigma^{\prime}$ is in ${\rm lk}_X(y)$. Since $\sigma \subset
\sigma^{\prime} \cup \{y\}$ and $y\in \sigma$, it follows that
$\sigma\in X$. \hfill $\Box$

\begin{theo}$\!\!\!${\bf .} \label{SC=SS}
Let $X$ be a normal pseudomanifold of dimension $d\geq 4$. Then
$X$ is a stacked sphere if and only if $X\in {\cal K}(d)$ and $X$
is simply connected.
\end{theo}

\noindent {\bf Proof.} If $X$ is a stacked sphere of dimension
$d\geq 2$ then $X$ is simply connected and $X \in {\cal K}(d)$ by
Lemma \ref{LSSSB5}. Conversely, let $X\in {\cal K}(d)$ be simply
connected and $d\geq 4$. We prove that $X$ is a stacked sphere by
induction on the number $n$ of vertices of $X$. If $n = d+2$ then
$X = S^{\,d}_{d+2}$ is a stacked sphere. So, assume $n > d+2$,
and we have the result for all smaller values of $n$. Now, take
an induced standard $(d-1)$-sphere $S$ in $X$ (Lemma
\ref{LSSSB10}). Let $\widetilde{X}$ be obtained from $X$ by
deleting the handle over $S$ (Lemma \ref{LEHD}). Clearly, since
$X$ is simply connected, $\widetilde{X}$ must be disconnected. If
$X_1$, $X_2$ are the connected components of $\widetilde{X}$,
then we have $X = X_1 \#X_2$. Clearly, $X_1$, $X_2$ are also
simply connected. Also, by Lemma \ref{LSSSB9} (b), $X_1, X_2 \in
{\cal K}(d)$. Hence by the induction hypothesis, $X_1$, $X_2$ are
stacked spheres. Therefore, by Lemma \ref{LSSSB8}, $X$ is a
stacked sphere. \hfill $\Box$

\bigskip

We shall not use this theorem in what follows. It is included
only for completeness.

\section{Generalized bistellar moves (GBMs)}

\begin{defn}$\!\!\!${\bf .} \label{DGBM1}
{\rm Let $X$ be a $d$-dimensional weak pseudomanifold. Let $B_1$,
$B_2$ be two combinatorial $d$-balls such that $B_1$ is a
subcomplex of $X$ and $\partial B_1 = \partial B_2 = B_2 \cap X$.
Then the pure $d$-dimensional simplicial complex $\widetilde{X} =
(X \setminus B_1) \cup B_2$ is said to be obtained from $X$ by a
{\em generalised bistellar move} (GBM) with respect to the pair
$(B_1, B_2)$. Observe that $\widetilde{X}$ is also a
$d$-dimensional weak pseudomanifold. [Let $\tau$ be a $(d -
1)$-face of $\widetilde{X}$. If $\tau \in B_2 \setminus
\partial B_2$ then $\tau$ is in two facets in $B_2$. If $\tau \in
\widetilde{X} \setminus B_2$ then $\tau$ is in two facets in $X
\setminus B_1 = \widetilde{X} \setminus B_2$. If $\tau \in
\partial B_1 =
\partial B_2$ then $\tau$ is in one facet in $X \setminus B_1 =
\widetilde{X} \setminus B_2$ and in one facet in $B_2$.] Notice
that we then have  $\partial B_2 = \partial B_1 = B_1 \cap
\widetilde{X}$, and $X$ is obtained from $\widetilde{X}$ by the
(reverse) generalised bistellar move with respect to the pair
$(B_2, B_1)$.  In case both $B_1$ and $B_2$ are $d$-balls with at
most $d + 2$ vertices (and hence at least one has $d+2$ vertices)
then this construction reduces to the usual bistellar move.
Clearly, if $\widetilde{X}$ is obtained from $X$ by a generalised
bistellar move then $|\widetilde{X}|$ is homeomorphic to $|X|$
and if the dimension of $X$ is at most 3 then $|\widetilde{X}|$
is pl homeomorphic to $|X|$.

}
\end{defn}

\begin{lemma}$\!\!\!${\bf .} \label{LGBM1}
If $\widetilde{X}$ is obtained from $X$ by a GBM, then
$\widetilde{X}$ is a normal pseudomanifold if and only if $X$ is
a normal pseudomanifold.
\end{lemma}

\noindent {\bf Proof.} Let $X$ be a normal pseudomanifold. We
prove that $\widetilde{X}$ is a normal pseudomanifold by induction
on the dimension $d$ of $X$. If $d = 1$ then the result is
trivial. Assume that the result is true for all normal
pseudomanifolds of dimension $< d$ and $X$ is a normal
pseudomanifold of dimension $d \geq 2$. Let $\widetilde{X}$ be
obtained from $X$ by a GBM with respect to the pair $(B_1, B_2)$.
Since $X$ is connected, it follows that $\widetilde{X}$ is
connected. We have observed that $\widetilde{X}$ is a weak
pseudomanifold. Let $\alpha$ be a face of dimension at most $d -
2$. If $\alpha \in B_2 \setminus
\partial B_2$ then ${\rm lk}_{\widetilde{X}}(\alpha) = {\rm
lk}_{B_2}(\alpha)$ is connected. If $\alpha \in \widetilde{X}
\setminus B_2$ then ${\rm lk}_{\widetilde{X}}(\alpha) =  {\rm
lk}_{X}(\alpha)$ is connected. If $\alpha \in \partial B_1 =
\partial B_2$ then ${\rm lk}_{\widetilde{X}}(\alpha)$ is obtained
from ${\rm lk}_{X}(\alpha)$ by the GBM with respect to the pair
$({\rm lk}_{B_1}(\alpha), {\rm lk}_{B_2}(\alpha))$. Since ${\rm
lk}_{X}(\alpha)$ is a normal pseudomanifold of dimension $< d$,
by induction hypothesis, ${\rm lk}_{\widetilde{X}}(\alpha)$ is a
normal pseudomanifold. In particular, ${\rm
lk}_{\widetilde{X}}(\alpha)$ is connected. This implies that
$\widetilde{X}$ is a normal pseudomanifold. Since $X$ is obtained
from $\widetilde{X}$ by the reverse GBM, the converse follows.
\hfill $\Box$

\begin{lemma}$\!\!\!${\bf .} \label{LGBM2}
Let $X$ be an $n$-vertex connected oriented triangulated
$2$-manifold. Then one of the following four cases must arise\,:
$(i)$ $X = S^{\,2}_4$,  $(ii)$ $X = X_1 \# X_2$ where $X_1$,
$X_2$ are connected orientable triangulated $2$-manifolds, $(iii)$
$X$ is obtained from a connected orientable triangulated
$2$-manifold $Y$ by an elementary handle addition or $(iv)$ for
each $u \in V(X)$, there exists a ball $B_u$ with $V(B_u) =
V({\rm lk}_X(u))$, $\partial B_u = B_u \cap X = {\rm lk}_X(u)$ so
that $X$ is obtained from the $(n - 1)$-vertex connected
orientable triangulated $2$-manifold $Y := (X\setminus {\rm
star}_X(u)) \cup B_u$  by the GBM with respect to the pair $(B_u,
{\rm star}_X(u))$.
\end{lemma}

\noindent {\bf Proof.} Assume that $X \neq S^{\,2}_4$. Take a
vertex $x$ of $X$. If ${\rm lk}_X(x)$ has a diagonal $yz$ which
is an edge of $X$, then the set $\{x, y, z\}$ induces an
$S^{\,1}_3$ in $X$. Since $X$ is orientable, this $S^{\,1}_3$ is
two sided. Let $Y$ be obtained from $X$ by a handle deletion over
this $S^{\,1}_3$ ($Y$ exists by Lemma \ref{LEHD}). Clearly, $Y$ is
also orientable. If $Y$ is connected then we are in Case $(iii)$
of this lemma. Otherwise, by Lemma \ref{LEHD}, $X = Y_1 \# Y_2$,
where $Y_1$, $Y_2$ are the connected components of $Y$. Here we
are in Case $(ii)$ of the lemma.

Finally, assume that none of the diagonals of the cycle ${\rm
lk}_X(x)$ are edges of $X$ for each $x\in V(X)$. Then, for each
$x\in V(X)$, $X$ is obtained from an $(n - 1)$-vertex triangulated
$2$-manifold $Y$ by a GBM with respect to $(B_x, \, {\rm
star}_X(x))$, where $B_x$ is any 2-ball with $V(B_x) = V({\rm
lk}_X(x))$ and $\partial B_x = {\rm lk}_X(x)$. Then we are in the
Case $(iv)$ of the lemma. \hfill $\Box$

\begin{remark}$\!\!\!${\bf .} \label{mintri}
{\rm Lemma \ref{LGBM2} shows, in particular, that any minimal
triangulation of a connected, orientable $2$-manifold of positive
genus must arise as the connected sum of two triangulated
$2$-manifolds or by handle addition over a triangulated
$2$-manifold of smaller genus. This fact should be useful in the
explicit classification of minimal triangulations of orientable
$2$-manifolds of small genus. Lemma \ref{LGBM2} also shows that
any triangulated $2$-sphere on $n$ ($ > 4$) vertices arises from
an $(n - 1)$-vertex triangulated $2$-sphere by a GBM. This should
help in simplifying the existing classifications and obtaining
new classifications of triangulated $2$-spheres with few
vertices. }
\end{remark}

\section{Gromov's combinatorial notion of rigidity}

Throughout this section, we use the following definition due to
Gromov (except that Gromov does not include connectedness as a
requirement for rigidity; but it seems anathema to call a
disconnected object rigid!). Thus $q$-rigidity hitherto refers to
Gromov's $q$-rigidity, without further mention.

\begin{defn}$\!\!\!${\bf .} \label{rigid}
{\rm Let $X$ be a $d$-dimensional simplicial complex and $q$ be a
positive integer. We shall say that $X$ is {\em $q$-rigid} if $X$
is connected and, for any set $A\subseteq V(X)$ which is disjoint
from at least one $d$-face of $X$, the number of edges of $X$
intersecting $A$ is $\geq mq$, where $m= \#(A)$. }
\end{defn}

\begin{lemma}$\!\!\!${\bf .} \label{LGR1}
Let $X$ be an $n$-vertex $d$-dimensional simplicial complex. If
$X$ is $q$-rigid then the number of edges of $X$ is $\geq (n - d -
1)q + {d + 1 \choose 2}$.
\end{lemma}

\noindent {\bf Proof.} Let $e$ be the number of edges of $X$. Fix
a $d$-face $\sigma$ of $X$ and put $A = V(X) \setminus \sigma$.
Then $\#(A) = n-d-1$ and exactly $e-{d+1\choose 2}$ edges
intersect $A$. \hfill $\Box$

\begin{defn}$\!\!\!${\bf .} \label{mrigid}
{\rm Let $X$ be an $n$-vertex $d$-dimensional simplicial complex
and $q$ a positive integer. We shall say that $X$ is {\em
minimally $q$-rigid} if $X$ is $q$-rigid and has exactly $(n - d
- 1)q + {d + 1 \choose 2}$ edges (i.e., if the lower bound in
Lemma \ref{rigid} is attained by $X$).}
\end{defn}

\begin{lemma}$\!\!\!${\bf .} \label{LGR2}
A connected simplicial complex is $q$-rigid if and only if the
cone over it is $(q+1)$-rigid. It is minimally $q$-rigid if and
only if the cone over it is minimally $(q+1)$-rigid.
\end{lemma}

\noindent {\bf Proof.} Let $X$ be an $n$-vertex $d$-dimensional
simplicial complex and $C(X) = x\ast X$ be the cone over $X$ with
cone-vertex $x$. Note that all the $(d +1)$-faces of $C(X)$ pass
through $x$, so that $A \subseteq V(C(X))$ is disjoint from a
$(d+1)$-face if and only if $A\subseteq V(X)$ and $A$ is disjoint
from a $d$-face of $X$. Also $C(X)$ has exactly $m = \#(A)$ more
edges than $X$ which intersect $A$ (viz., the edges joining $x$
with the vertices of $A$). In consequence, the number of edges of
$X$ intersecting $A$ is $\geq mq$ if and only if the number of
edges of $C(X)$ intersecting $A$ is $\geq m(q + 1)$. This proves
the first part. The second part follows since $C(X)$ has one more
vertex and $n$ more edges than $X$. \hfill $\Box$

\begin{lemma}$\!\!\!${\bf .} \label{LGR3}
Let $X_1$, $X_2$ be subcomplexes of a simplicial complex $X$ such
that $X = X_1 \cup X_2$ and $\dim(X_1\cap X_2) = \dim(X)$. If
$X_1$, $X_2$ are both $q$-rigid then $X$ is $q$-rigid. If,
further, $X$ is minimally $q$ rigid then both $X_1$, $X_2$ are
minimally $q$-rigid.
\end{lemma}

\noindent {\bf Proof.} Since $X_1$, $X_2$ are both connected, our
assumption implies that $X$ is connected. Let $\dim(X)= d$. Since
$\dim(X_1\cap X_2) = \dim(X)$, it follows that $\dim(X_1) =
\dim(X_2) = \dim(X_1\cap X_2) = d$. Let $A\subseteq V(X)$ be
disjoint from some $d$-face $\sigma \in X = X_1\cup X_2$. Without
loss of generality, $\sigma \in X_1$. Write $A_1 = A\cap V(X_1)$
and $A_2 = A\setminus V(X_1)$. Say $m = \#(A)$, $m_i = \#(A_i)$,
$i = 1, 2$. Thus, $m = m_1 +m_2$. Note that $A_1 \subseteq
V(X_1)$ is disjoint from the $d$-face $\sigma$ of $X_1$. Also, if
$\tau$ is a $d$-face of $X_1\cap X_2$, then $\tau$ is a $d$-face
of $X_2$ disjoint from $A_2$ (since $\tau\subseteq V(X_1)$ and
$A_2$ is disjoint from $V(X_1)$). Since, $X_1$, $X_2$ are
$q$-rigid, we have at least $m_1q$ edges of $X_1$ meeting $A_1$
and at least $m_2q$ edges of $X_2$ meeting $A_2$. Also, as
$V(X_1)$ and $A_2$ are disjoint, no edge of $X_1$ meets $A_2$.
Therefore, we have at least $m_1q+m_2q = mq$ distinct edges of
$X$ meeting $A$. This proves that $X$ is $q$-rigid.

Now, if $X$ is minimally $q$-rigid, then taking $A$ to be the
complement in $V(X)$ of a $d$-face of $X_1$, one gets exactly
$mq$ edges of $X$ meeting $A$. Since we have equality in the
above argument, it follows that exactly $m_1q$ edges of $X_1$
intersect $A_1 = A\cap V(X_1)$. Since $A_1$ is the complement in
$V(X_1)$ of a $d$-face of $X_1$, this shows that $X_1$ is then
minimally $q$-rigid. Since the assumptions are symmetric in $X_1$
and $X_2$, in this case $X_2$ is also minimally $q$-rigid. \hfill
$\Box$

\begin{lemma}$\!\!\!${\bf .} \label{LGR4}
Let $\{X_{\alpha} : \alpha \in I\}$ be a finite family of
$q$-rigid subcomplexes of a simplicial complex $X$. Suppose there
is a connected graph $H$ with vertex set $I$ such that whenever
$\alpha, \beta \in I$ are adjacent in $H$, we have
$\dim(X_{\alpha} \cap X_{\beta}) =\dim(X)$. Also suppose
$\cup_{\alpha\in I}X_{\alpha} = X$. Then $X$ is $q$-rigid. If,
further, $X$ is minimally $q$-rigid, then each $X_{\alpha}$ is
minimally $q$-rigid.
\end{lemma}

\noindent {\bf Proof.} Induction on $\#(I)$. If $\#(I) = 1$ then
the result is trivial. For $\#(I)= 2$, the result is just Lemma
\ref{LGR3}. So suppose $\#(I) >2$ and we have the result for
smaller values of $\#(I)$. Since $H$ is a connected graph, there
is $\alpha_0\in I$ such that the induced subgraph of $H$ on the
vertex set $I\setminus\{\alpha_0\}$ is connected (for instance,
one may take $\alpha_0$ to be an end vertex of a spanning tree in
$H$). Applying the induction hypothesis to the family
$\{X_{\alpha} : \alpha \neq \alpha_0\}$, one gets that $Y_1 =
\cup_{\alpha \neq \alpha_0} X_{\alpha}$ is $q$-rigid. Since $Y_2
= X_{\alpha_0}$ is also $q$-rigid, $X = Y_1\cup Y_2$, and
$\dim(Y_1\cap Y_2) = \dim(X)$ (if $\alpha_0$ is adjacent to
$\alpha_1$ in $H$ then $\dim(X) \geq \dim(Y_1\cap Y_2) \geq
\dim(X_{\alpha_1}\cap Y_2) = \dim(X)$), induction hypothesis (or
Lemma \ref{LGR3}) implies that $X$ is $q$-rigid. Now, if $X$ is
minimally $q$-rigid then, by Lemma \ref{LGR3}, so are $Y_1$ and
$Y_2$. Since $Y_1$ is minimally $q$-rigid, induction hypothesis
then implies that $X_{\alpha}$ is minimally $q$-rigid for
$\alpha\neq \alpha_{0}$ (and also for $\alpha = \alpha_0$ since
$X_{\alpha_0} = Y_2$). \hfill $\Box$

\begin{lemma}$\!\!\!${\bf .} \label{LGR5}
Let $X$ be a connected pure $d$-dimensional simplicial complex.
$(a)$ If each vertex link of $X$ is $q$-rigid then $X$ is $(q +
1)$-rigid. $(b)$ If, further, $X$ is minimally $(q+1)$-rigid then
all the vertex links of $X$ are minimally $q$-rigid.
\end{lemma}

\noindent {\bf Proof.} Let $I = V(X)$ and $H$ be the edge graph
of $X$. Since $X$ is connected, so is $H$. For $\alpha\in I$,
${\rm st}(\alpha)$ is a cone over the $q$-rigid complex ${\rm
lk}(\alpha)$, and hence by Lemma \ref{LGR2}, ${\rm st}(\alpha)$
is $(q+1)$-rigid for each $\alpha\in I$. Since $X$ is pure, the
family $\{{\rm st}(\alpha) : \alpha \in I\}$ satisfies the
hypothesis of Lemma \ref{LGR4}. Hence $X$ is $(q+1)$-rigid. If it
is minimally $(q+1)$-rigid, then by Lemma \ref{LGR4}, each ${\rm
st}(\alpha)$ is minimally $(q+1)$-rigid, and hence, by Lemma
\ref{LGR2}, ${\rm lk}(\alpha)$ is minimally $q$-rigid for all
$\alpha\in I$. \hfill $\Box$

\begin{lemma}$\!\!\!${\bf .} \label{LGR6}
Let $X_1$, $X_2$ be $d$-dimensional normal pseudomanifolds. If
$X_1$, $X_2$ are $(d + 1)$-rigid then their elementary connected
sum $X_1 \# X_2$ is $(d + 1)$-rigid. If, further, $X_1 \# X_2$ is
minimally $(d + 1)$-rigid then both $X_1$ and $X_2$ are minimally
$(d + 1)$-rigid.
\end{lemma}

\noindent {\bf Proof.} Since $X_1$, $X_2$ are both connected, so
is $X_1 \# X_2$. Let $\sigma_i$ be a facet of $X_i$ ($i = 1, 2$)
and $f \colon \sigma_1 \to \sigma_2$ be a bijection, such that $X
= X_1 \# X_2$ is obtained from $X_1 \sqcup X_2 \setminus
\{\sigma_1, \sigma_2\}$ via an identification through $f$. We view
$V(X_i)$ as a subset of $V(X)$ in the obvious fashion. Put
$\widetilde{X} = (X_1 \# X_2) \cup \{\sigma_1 = \sigma_2\}$. Then
$X_1$, $X_2$ are subcomplexes of $\widetilde{X}$ satisfying the
hypothesis of Lemma \ref{LGR3} with $q = d+1$. Hence, by Lemma
\ref{LGR3}, $\widetilde{X}$ is $(d+1)$-rigid. Since $X_1 \# X_2$
is a subcomplex of $\widetilde{X}$ of the same dimension with the
same set of edges, it follows that $X_1 \# X_2$ is $(d+1)$-rigid.

If $X_1 \# X_2$ is minimally $(d+1)$-rigid, then so is
$\widetilde{X}$ and hence, by Lemma \ref{LGR3}, so are $X_1$,
$X_2$. \hfill $\Box$

\begin{lemma}$\!\!\!${\bf .} \label{LGR7}
Let $Y$ be a $d$-dimensional normal pseudomanifold which is
obtained from a $d$-dimensional normal pseudomanifold $X$ by an
elementary handle addition. If $X$ is $(d + 1)$-rigid then $Y$ is
$(d + 1)$-rigid.
\end{lemma}

\noindent {\bf Proof.} Let $Y = X^{\psi}$, where $\psi \colon
\sigma_1 \to \sigma_2$ is an admissible bijection between two
disjoint facets $\sigma_1$, $\sigma_2$ of $X$. Thus $Y$ is
obtained from $X\setminus \{\sigma_1, \sigma_2\}$ by identifying
$x$ with $\psi(x)$ for each $x\in \sigma_1$ (cf. Definition
\ref{DEHA}). Let's  identify $V(Y)$ with $V(X) \setminus \sigma_2$
via the quotient map $V(X) \to V(Y)$. Let $A \subseteq V(Y)$ be
an $m$-set disjoint from a facet $\sigma$ of $Y$. Then, under
this identification $A \subseteq V(X)$ is disjoint from $\sigma$
and it follows from the definition of $X^{\psi}$ that $\sigma$ is
a facet of $X$. This implies, by $(d+1)$-rigidity of $X$, that at
least $m(d+1)$ edges of $X$ meet $A$. Since $A\cap \sigma_2 =
\emptyset$, these edges corresponds to distinct edges of $Y$
under our identification. Hence $Y$ is $(d+1)$-rigid. \hfill
$\Box$

\begin{lemma}$\!\!\!${\bf .} \label{LGR8}
Let $X$ be a triangulated $2$-manifold. Suppose for each vertex
$u$ of $X$, there is a triangulated $2$-manifold $X_u$ with
vertex-set $V(X)\setminus \{u\}$, and a triangulated $2$-ball $B_u
\subseteq X_u$ with vertex-set $V({\rm lk}_X(u))$ such that $X$
is obtained from $X_u$ by the GBM with respect to the pair $(B_u,
{\rm star}_X(u))$. If $X_u$ is $3$-rigid for all $u \in V(X)$,
then $X$ is $3$-rigid.
\end{lemma}

\noindent {\bf Proof.} Take any set $A\subseteq V(X)$ which is
disjoint from at least one 2-face $\sigma$ of $X$. Say $\#(A) =
m$. Fix a vertex $x\in A$, say of degree $k$. Take a 2-ball $B$
with vertex set $V(B) = V({\rm lk}(x))$ as in the hypothesis.
Note that $B$ is a $k$-vertex 2-ball with $k$ edges in the
boundary (viz., the edges of ${\rm lk}_X(x)$), hence it has $k-3$
edges in the interior: these are not edges of $X$. By assumption
$X_x = (X \setminus {\rm st}(x)) \cup B$ is 3-rigid, so that at
least $3(m - 1)$ edges of $X_x$ intersect $\widetilde{A}$, and
hence also $A$. Of these edges, at most $k - 3$ edges are not in
$X$. Thus at least $3(m - 1) - (k - 3)$ edges of $X$ (not passing
through $x$) meet $A$. Also, all the $k$ edges of $X$ through $x$
meet $A$. Thus we have a total of at least $3(m-1) - (k-3) +k =
3m$ edges of $X$ meeting $A$. Hence $X$ is 3-rigid. \hfill $\Box$

\section{$(d+1)$-rigidity of normal $d$-pseudomanifolds}

\begin{lemma}$\!\!\!${\bf .} \label{LRNP1}
Let $X$ be a $2$-dimensional normal pseudomanifold. Then $X$ is
$3$-rigid. $X$ is minimally $3$-rigid if and only if $X$ is a
triangulated $2$-sphere.
\end{lemma}

\noindent {\bf Proof.} Since $X$ is 2-dimensional normal
pseudomanifold, it follows that $X$ is a connected triangulated
$2$-manifold.

First assume that $X$ is orientable. Recall that the connected
orientable closed 2-manifolds are classified up to homeomorphism
by their genus $g$. The genus is related to the Euler
characteristic $\chi$ by the formula $\chi = 2 - 2g$. With any
$X$ as above, we associate the parameter $(g, n)$, where $g$ is
the genus of $|X|$ and $n$ is the number of vertices of $X$.
Let's well order the collection of all possible parameters by the
lexicographic order $\prec$. That is, $(g_1, n_1) \prec (g_2,
n_2)$ if either $g_1 < g_2$ or else $g_1 = g_2$ and $n_1 < n_2$.
We prove the 3-rigidity of $X$ by induction with respect to
$\prec$. Notice that the smallest parameter is $(0, 4)$
corresponding to $X = S^{\,2}_4$, which is trivially 3-rigid.
This starts the induction. If $(g, n) \succ (0, 4)$, then $X$ is
as in Case $(ii)$, $(iii)$ or $(iv)$ of Lemma \ref{LGBM2}.

If $X$ is as in $(ii)$, then $X = X_1 \# X_2$ where $X_1$, $X_2$
are connected orientable $2$-manifold with small parameters. Hence
by induction hypothesis, $X_1$, $X_2$ are 3-rigid. Hence by Lemma
\ref{LGR6}, $X$ is 3-rigid. If $X$ is as in Case $(iii)$, then
$X$ is obtained from a connected orientable triangulated
2-manifold $Y$ of smaller genus, by elementary handle addition.
By induction hypothesis, $Y$ is 3-rigid, and hence by Lemma
\ref{LGR7}, $X$ is 3-rigid. If $X$ is as in Case $(iv)$ of Lemma
\ref{LGBM2}, then it satisfies the hypothesis of Lemma
\ref{LGR8}, and hence is 3-rigid. This completes the induction.

Now suppose $X$ is non-orientable. Let $\widehat{X}$ be the
orientable double cover of $X$. By the above, $\widehat{X}$ is
3-rigid. Since the covering map $V(\widehat{X}) \to V(X)$ is a
two-to-one simplicial map, it is immediate that $X$ is 3-rigid.

Finally, $X$ is minimally 3-rigid $\Longleftrightarrow$ number of
edges in $X$ is $3(n-2)$ $\Longleftrightarrow$ the Euler
characteristic of $X$ is 2 $\Longleftrightarrow$ $X$ is a
triangulated 2-sphere. \hfill $\Box$

\begin{prop}$\!\!\!${\bf .} \label{P7.1}
Let $X$ be a $d$-dimensional normal pseudomanifold. If $d \geq 2$
then $X$ is $(d + 1)$-rigid. If, further, $d \geq 3$ and $X$ is
minimally $(d + 1)$-rigid, then all the vertex links of $X$ are
minimally $d$-rigid.
\end{prop}

\noindent {\bf Proof.} The proof is by induction on $d$. For $d =
2$ this is Lemma \ref{LRNP1}. For $d \geq 3$, all the vertex
links of $X$ are $(d - 1)$-dimensional normal pseudomanifolds and
hence, by the induction hypothesis, all vertex links of $X$ are
$d$-rigid. So the result follows from Lemma \ref{LGR5}. \hfill
$\Box$

\begin{lemma}$\!\!\!${\bf .} \label{LRNP2}
Let $X$ be a minimally $(d + 1)$-rigid normal pseudomanifold of
dimension $d \geq 3$. Then every clique of size $\leq d$ in the
edge graph of $X$ is a face of $X$.
\end{lemma}

\noindent {\bf Proof.} Let $I = V(X)$ and let $H$ be the edge
graph of $X$. For $\alpha\in I$, let $H_{\alpha}$ be the induced
subgraph of $H$ on the vertex-set $V({\rm lk}(\alpha))$ and put
$X_{\alpha} = {\rm st}(\alpha)\cup H_{\alpha}$. By Lemma
\ref{LGR2} and Theorem \ref{P7.1}, ${\rm st}(\alpha)$ is
$(d+1)$-rigid and hence so is $X_{\alpha}$. Thus $\{X_{\alpha} :
\alpha\in I\}$ satisfies the hypothesis of Lemma \ref{LGR4}.
Since $X$ is minimally $(d+1)$-rigid, it follows that
$X_{\alpha}$ is minimally $(d+1)$-rigid for each $\alpha\in I$.
But $X_{\alpha}\supseteq {\rm st}(\alpha)$, $V(X_{\alpha}) =
V({\rm st}(\alpha))$ and ${\rm st}(\alpha)$ is $(d+1)$-rigid.
Therefore, $X_{\alpha}$ and ${\rm st}(\alpha)$ have the same edge
graph. That is, $H_{\alpha} \subseteq {\rm st}(\alpha)$. Thus,
each clique of size $\leq 3$ through $\alpha$ is a face of $X$.
Since this holds for each $\alpha\in I$, it follows that each
clique of size $\leq 3$ in $H$ is a face of $X$.

Now, by an induction on $k$, one sees that for $k\leq d$, any
$k$-clique of $H$ is a face of $X$: if $C$ is a $k$-clique (and
$k\geq 4$ and hence $d\geq 4$), then for any $x\in C$, $C
\setminus \{x\}$ is a $(k-1)$-clique of ${\rm lk}(x)$ and
$\dim({\rm lk}(x)) = d-1\geq 3$. Therefore, $C \setminus \{x\}$ is
a face of ${\rm lk}(x)$ and hence $C$ is a face of $X$. \hfill
$\Box$

\begin{lemma}$\!\!\!${\bf .} \label{LRNP3}
Let $X$ be a minimally $(d+1)$-rigid normal pseudomanifold of
dimension $d \geq 3$. Then the edge graph of $X$ has a clique of
size $d + 2$.
\end{lemma}

\noindent {\bf Proof.} If we have the result for $d=3$ then the
result follows for all $d\geq 3$ by a trivial induction on
dimension (using the second statement in Proposition \ref{P7.1}).
So, we may assume $d=3$.

Let $n \geq 5$ be the number of vertices of $X$. Since $X$ is
minimally 4-rigid, it has $4n - 10$ edges and hence the average
degree of the vertices is $\frac{2(4n - 10)}{n} < 8$. Therefore,
$X$ has a vertex $x$ of degree $\leq 7$. Then, by Lemmas
\ref{LGR5} and \ref{LRNP1}, ${\rm lk}(x)$ is a triangulated
2-sphere on $\leq 7$ vertices. If possible, suppose ${\rm lk}(x)$
has no vertex of degree 3. It is easy to see that up to
isomorphism there are only two such $S^{\hspace{.2mm}2}$, namely
$S^{\hspace{.2mm}0}_2 \ast S^{\hspace{.2mm}1}_m$ with $m = 4$ or
5. Thus ${\rm lk}(x)$ is one of these two spheres, say ${\rm
lk}(x) = S^{\hspace{.2mm}0}_2(\{y, z\}) \ast
S^{\hspace{.2mm}1}_m(A)$. Since $xyz$ is not a 2-face, by Lemma
\ref{LRNP2}, $yz$ is not an edge of $X$. Put $B_1 = {\rm
st}_X(x)$, $B_2 = B^1_2(\{x, y\}) \ast S^{\hspace{.2mm}1}_m(A)$.
Set $\widetilde{X} = (X \setminus B_1) \cup B_2$. Then
$\widetilde{X}$ is obtained from $X$ by a GBM. Hence
$\widetilde{X}$ is a 3-dimensional normal pseudomanifold with $n
- 1$ vertices and $4n - 10 - (m + 2) + 1 = 4n - 11- m < 4(n - 1)
- 10$ edges (as $m \geq 4$). This is impossible since
$\widetilde{X}$ is 4-rigid by Proposition \ref{P7.1}. This proves
that ${\rm lk}(x)$ has a vertex $y$ of degree 3. Then the
vertex-set of ${\rm st}(xy)$ is a 5-clique. This completes the
proof. \hfill $\Box$

\begin{lemma}$\!\!\!${\bf .} \label{LRNP4}
Let $X$ be an $n$-vertex minimally $(d + 1)$-rigid
$d$-dimensional normal pseudomanifold. If $d \geq 3$ and $n
> d + 2$ then $X$ contains a standard $(d - 1)$-sphere $S$ as an
induced subcomplex.
\end{lemma}

\noindent {\bf Proof.} By Lemma \ref{LRNP3}, there is a $(d +
2)$-set $C \subseteq V(X)$ which is a clique of the edge graph of
$X$. If all the $(d + 1)$-subsets of $C$ were facets of $X$ then
the induced subcomplex of $X$ on the vertex-set $C$ would be a
proper subcomplex which is a (standard) $d$-sphere. This is not
possible since $X$ is a $d$-dimensional normal pseudomanifold. So,
there is a $(d + 1)$-set $C_0 \subseteq C$ such that $C_0$ is not
a facet of $X$. But $C_0$ is a $(d + 1)$-clique of the edge
graph  of $X$, so by Lemma \ref{LRNP2}, all proper subsets of
$C_0$ are faces of $X$. Thus the induced subcomplex $S$ of $X$ on
the vertex-set $C_0$ is a standard $(d-1)$-sphere. \hfill $\Box$

\begin{lemma}$\!\!\!${\bf .} \label{LRNP5}
If $X$ is a minimally $4$-rigid $3$-dimensional normal
pseudomanifold then $X$ is a stacked $3$-sphere.
\end{lemma}

\noindent {\bf Proof.} By Theorem \ref{P7.1}, all the vertex
links are minimally $3$-rigid. Therefore, by Lemma \ref{LRNP1},
$X$ is a triangulated 3-manifold. Let the number of vertices in
$X$ be $n$. We wish to prove by induction on $n$ that $X$ must be
a stacked 3-sphere. This is trivial for $n = 5$, so that we may
assume that $n > 5$ and we have the result for smaller values of
$n$.

By Lemma \ref{LRNP4}, $X$ contains a standard $2$-sphere $S$ as
an induced subcomplex. Since $S$ is a 2-sphere, $S$ is two-sided
in $X$. Let $Y$ be the simplicial complex obtained from $X$ by
deleting the ``handle'' over $S$. Since $X$ is a triangulated
$3$-manifold, by Lemma \ref{LSSSB9} $(a)$, $Y$ is a triangulated
$3$-manifold. Also, $Y$ has $n + 4$ vertices and $4n - 10 + {4
\choose 2} < 4(n + 4) - {5 \choose 2}$ edges. Therefore $Y$ is
not $4$-rigid and hence, by Theorem \ref{P7.1}, $Y$ must be
disconnected. Since $X$ is connected, Lemma \ref{LEHD} implies
that  $X = Y_1 \# Y_2$, where $Y_1$, $Y_2$ are $3$-dimensional
normal pseudomanifolds. Since $X$ is minimally $4$-rigid, Lemma
\ref{LGR6} implies that $Y_1$, $Y_2$ are both minimally
$4$-rigid. Let $Y_i$ have $n_i$ vertices ($i = 1, 2$). Since $n_1
+ n_2 = n + 4$, $n_1 > 4$, $n_2 > 4$, it follows that $n_1 < n$,
$n_2 < n$. Therefore, by induction hypothesis, $Y_1$, $Y_2$ are
stacked $3$-spheres. Since $X$ is an elementary connected sum of
$Y_1$ and $Y_2$, Lemma \ref{LSSSB8} $(b)$ implies that $X$ is a
stacked $3$-sphere. \hfill $\Box$

\begin{prop}$\!\!\!${\bf .} \label{P7.2}
For $d \geq 3$, the stacked $d$-spheres are the only minimally
$(d + 1)$-rigid $d$-dimensional normal pseudomanifolds.
\end{prop}

\noindent {\bf Proof.} If $X$ is an $n$-vertex stacked $d$-sphere
then (cf. Lemma \ref{LSSSB2}) the number of edges of $X$ is $(d +
1)n - {d + 2 \choose 2}$, so that $X$ is minimally $(d +
1)$-rigid by Theorem \ref{P7.1}.

For the converse, let $X$ be a minimally $(d + 1)$-rigid
$d$-dimensional normal pseudomanifold, with $d \geq 3$. We prove
by induction on $d$ that $X$ is a stacked $d$-sphere. The $d = 3$
case is Lemma \ref{LRNP5}. So, assume $d > 3$ and we have the
result for smaller values of $d$. By Theorem \ref{P7.1} and
induction hypothesis, all the vertex links of $X$ are stacked $(d
- 1)$-spheres. That is, $X$ is in the class ${\cal K}(d)$ (cf.
Definition \ref{DSS2}). In particular, $X$ is a triangulated
$d$-manifold.

Let the number of vertices in $X$ be $n$. We wish to prove by
induction on $n$ that $X$ must be a stacked $d$-sphere. This is
trivial for $n = d + 2$, so that we may assume that $n > d + 2$
and we have the result for smaller values of $n$.

By Lemma \ref{LRNP4} (also by Lemma \ref{LSSSB10}), $X$ contains a
standard $(d-1)$-sphere $S$ as an induced subcomplex. Since $d >
3$, $S$ is two-sided in $X$. Let $Y$ be the simplicial complex
obtained from $X$ by deleting the ``handle'' over $S$. Since $X$
is in the class ${\cal K}(d)$, by Lemma \ref{LSSSB9} $(b)$, $Y$
is in the class ${\cal K}(d)$. In particular, $Y$ is a
triangulated $d$-manifold. Also, $Y$ has $n + d + 1$ vertices and
$((d + 1)n - {d + 2 \choose 2}) + {d + 1 \choose 2} = (n + d +
1)(d + 1) - (d + 1)(d + 2) < (n + d + 1)(d + 1) - {d + 2 \choose
2}$ edges. Therefore $Y$ is not $(d + 1)$-rigid and hence, by
Theorem \ref{P7.1}, $Y$ must be disconnected. Since $X$ is
connected, Lemma \ref{LEHD} implies that  $X = Y_1 \# Y_2$, where
$Y_1$, $Y_2$ are $d$-dimensional normal pseudomanifolds. Since
$X$ is minimally $(d + 1)$-rigid, Lemma \ref{LGR6} implies that
$Y_1$, $Y_2$ are both minimally $(d + 1)$-rigid. Let $Y_i$ have
$n_i$ vertices ($i = 1, 2$). Since $n_1 + n_2 = n + d + 1$, $n_1
> d + 1$, $n_2 > d + 1$, it follows that $n_1 < n$, $n_2 < n$.
Therefore, by induction hypothesis, $Y_1$, $Y_2$ are stacked
$d$-spheres. Since $X$ is an elementary connected sum of $Y_1$
and $Y_2$, Lemma \ref{LSSSB8} $(b)$ implies that $X$ is a stacked
$d$-sphere.  \hfill $\Box$

\begin{theo}$\!\!\!${\bf .} \label{GRT}
For $d \geq 2$, all $d$-dimensional normal $d$-pseudomanifold are
$(d + 1)$-rigid. For $d \geq 3$, the stacked $d$-spheres are the
only minimally $(d + 1)$-rigid $d$-dimensional normal
pseudomanifolds.
\end{theo}

\noindent {\bf Proof.} Immediate from Propositions \ref{P7.1} and
\ref{P7.2}. \hfill $\Box$

\section{LBT for normal pseudomanifolds}

Now we are ready to state and prove the main result of this paper:

\begin{theo}$\!\!\!${\bf .} \label{LBT}
Let $X$ be any $d$-dimensional normal pseudomanifold. Then the
face-vector of $X$ satisfies;
$$
f_j(X) \geq \left\{\begin{array}{ll}
{d + 1 \choose j}f_0(X) - j{d + 2 \choose j + 1}, & \mbox{if }
~~ 1 \leq j < d, \\[1.5mm]
df_0(X) - (d+2)(d-1), & \mbox{if } ~~ j=d.
\end{array}
\right.
$$
Further, for $d \geq 3$, equality holds here for some $j$ if and
only if $X$ is a stacked sphere.
\end{theo}

\noindent {\bf Proof.} This is trivial for $d=1$. So, assume
$d>1$. For $j=1$, the result is immediate from Lemma \ref{LGR1},
Definition \ref{mrigid} and Theorem \ref{GRT}. So let $1 < j\leq
d$. Counting in two ways the incidences between vertices and
$j$-faces of $X$, we obtain
$$
f_j(X) = \frac{1}{j+1}\sum_{v\in V(X)} f_{j-1}({\rm lk}_X(v)).
$$
Since ${\rm lk}_X(v)$ is a $(d-1)$-dimensional normal
pseudomanifold with $\deg(v)$ vertices, induction hypothesis (on
the dimension) implies that
$$
f_{j-1}({\rm lk}_X(v)) \geq \left\{\begin{array}{ll} {d \choose j
- 1}\deg(v) - (j-1){d + 1 \choose j}, & \mbox{if }
~~ 1 < j < d, \\[1.5mm]
(d - 1)\deg(v) - (d+1)(d-2), & \mbox{if } ~~ j=d.
\end{array}
\right.
$$
Adding this inequality over all vertices $v$, and noting that
$\sum_{v\in V(X)} \deg(v) = 2 f_1(X)$, we conclude\,:
$$
f_j(X) \geq \left\{\begin{array}{ll}
 \frac{1}{j+1}\left(2{d \choose j-1}f_1(X) - (j-1){d + 1 \choose
 j}f_0(X)\right), & \mbox{if } ~~ 1 < j < d, \\[1.5mm]
\frac{1}{d+1}\left(2(d-1)f_1(X) - (d+1)(d-2)f_0(X)\right), &
\mbox{if } ~~ j=d.
\end{array}
\right.
$$
But, by the $j=1$ case of the theorem, $f_1(X) \geq (d+1)f_0(X) -
{d+2 \choose 2}$. Hence we get\,:
$$
f_j(X) \geq \left\{\begin{array}{ll}
 \frac{1}{j+1}\left(\left(2{d \choose j-1}(d+1) - (j-1){d + 1 \choose
 j}\right)f_0(X) -2{d \choose j-1}{d+2\choose 2}\right), & \mbox{if }  1 < j < d, \\[1.5mm]
\frac{1}{d+1}\left((2(d-1)(d+1) - (d+1)(d-2))f_0(X) - 2(d-1){d+2
\choose 2}\right), & \mbox{if }  j=d.
\end{array}
\right.
$$
Since $(d+1){d\choose j-1} = j{d+1\choose j}$ and ${d \choose
j-1}{d+2 \choose 2} = {d+2 \choose j+1}{j+1 \choose 2}$, this
inequality simplifies to the one stated in the theorem. From this
argument, it is clear that if the equality holds for some $j$,
then it also holds with $j=1$, so that (when $d \geq 3$) $X$ is a
stacked sphere in the case of equality. The converse is immediate
from Lemma \ref{LSSSB2}. \hfill $\Box$

\begin{remark}$\!\!\!${\bf .} 
{\rm The argument in the above proof (reducing the inequality for
arbitrary $j$ to the case $j=1$) is known as the M-P-W reduction
- after its independent inventors McMullen, Perles and Walkup.}
\end{remark}

\section{Some more lower bound conjectures}

\begin{defn}$\!\!\!${\bf .} \label{k-stacked}
{\rm For $0\leq k\leq d$, a triangulated $d$-sphere $X$ is said
to be a {\em $k$-stacked sphere} if there is a triangulated
$(d+1)$-ball $B$ such that $\partial B = X$ and ${\rm
skel}_{d-k}(B) = {\rm skel}_{d-k}(X)$. Recall that ${\rm
skel}_{d-k}(X)$, for instance, is the subcomplex of $X$
consisting of all its faces of dimension at most $d-k$.}
\end{defn}

\begin{defn}$\!\!\!${\bf .} \label{DGLBC2}
{\rm Let $X$ be a $d$-dimensional pseudomanifold and $u$ be a
vertex of $X$. Then, for a new symbol $v\notin V(X)$, the
$(d+1)$-dimensional pseudomanifold $\Sigma_{u, v}(X) : = (u\ast
{\rm ast}_X(u)) \cup (v\ast X)$ is called an {\em one point
suspension} of $X$. The geometric carrier of $\Sigma_{u, v}(X)$
is the suspension of $|X|$. In particular, $\Sigma_{u, v}(X)$ is
a triangulated $(d+1)$-sphere if $X$ is a triangulated $d$-sphere
(cf. \cite{bd2}).}
\end{defn}

\begin{lemma}$\!\!\!${\bf .} \label{LGLBC1}
If $X$ is a triangulated $d$-sphere then there is a triangulated
$(d+1)$-ball $\widetilde{X}$ such that $V(\widetilde{X}) = V(X)$
and $\partial \widetilde{X} = X$.
\end{lemma}

\noindent {\bf Proof.} Fix a vertex $u$ of $X$, and let $X_u = u
\ast {\rm ast}_X(u)$. Since $X$ is a triangulated $d$-sphere, it
follows that $\Sigma_{u, v}(X)$ is a triangulated $(d+1)$-sphere.
Thus, $X_u$ is the antistar of the vertex $v$ in the triangulated
$(d + 1)$-sphere $\Sigma_{u, v}(X)$ and the link of $v$ in
$\Sigma_{u, v}(X)$ is $X$. Therefore Lemma \ref{LSSSB1} implies
that $X_u$ is a triangulated $(d + 1)$-ball. Clearly, $V(X_u) =
V(X)$ and $\partial X_u = X$. Thus $\widetilde{X} = X_u$ works
for any vertex $u$ of $X$. \hfill $\Box$

\begin{remark}$\!\!\!${\bf .}
{\rm Trivially, for $0\leq k < l \leq d$, every $k$-stacked
$d$-sphere is also $l$-stacked. Further, the standard sphere
$S^{\,d}_{d+2}$ is the only 0-stacked $d$-sphere, while Lemma
\ref{LGLBC1} shows that all triangulated $d$-spheres are
$d$-stacked. Remark \ref{1-stacked} (a) shows that every stacked
sphere is 1-stacked. Conversely, the case $k=1$ of the following
proposition shows that the face-vector of any 1-stacked sphere
satisfies the LBT with equality, so that 1-stacked spheres are
precisely the stacked spheres. }
\end{remark}

\begin{prop}$\!\!\!${\bf .}
Let $k \geq 0$. Then for $d\geq 2k+1$, the $k$ components $f_0,
\dots, f_{k-1}$ of the face-vector of any $k$-stacked $d$-sphere
determines the rest of its face-vector by the formulae
$$
f_j = \left\{\begin{array}{l} {\displaystyle \sum_{i = -1}^{k - 1}
(-1)^{k-i+1}{j-i-1 \choose j-k}{d -i + 1  \choose j-i}f_{i}},
\hspace{3mm} \mbox{if } ~~ k \leq j \leq d-k, \\
{\displaystyle \sum_{i = -1}^{k - 1} (-1)^{k - i + 1} \left[
{j-i-1 \choose j-k}{d -i + 1  \choose j-i} -{k \choose d-j+1}
{d-i \choose d-k+1} \right.} \\
\mbox{} \hspace{27mm} + {\displaystyle \sum_{l = d-j}^{k - 1}
(-1)^{k - l} \left. {l \choose d-j}{d -i \choose
d-l+1}\right]f_{i}}, \hspace{3mm} \mbox{if } ~~ d-k+1 \leq j \leq
d.
\end{array}
\right.
$$
$($Here $f_{-1} = 1$, consistent with the convention that the
empty face is the only face of dimension $-1$ in any simplicial
complex.$)$
\end{prop}

\noindent {\bf Sketch of proof.} Let $X$ be a $k$-stacked
$d$-sphere. Let $B$ be a $(d+1)$-ball as in Definition
\ref{k-stacked}. Put $\widetilde{X} = B\cup (x \ast X)$, where $x$
is a new symbol. Thus $\widetilde{X}$ is a triangulated
$(d+1)$-sphere. Let $(f_0, f_1, \dots, f_d)$ and $(\tilde{f}_0,
\tilde{f}_1, \dots, \tilde{f}_{d+1})$ be the face-vectors of $X$
and $\widetilde{X}$, respectively. From the relation between $X$
and $B$, we get
\begin{eqnarray}
\tilde{f}_j = f_j + f_{j-1}, & 0\leq j\leq d-k.
\end{eqnarray}

Being triangulated spheres of dimension $d$ and $d+1$
respectively, $X$ and $\widetilde{X}$ satisfy the following
Dehn-Sommerville equations (cf. \cite[9.2.2, Page 148]{g})\,:
\begin{eqnarray}
\sum_{i=-1}^{j-1}(-1)^{d-i-1}{d-i \choose d-j+1}f_i  =
\sum_{i=-1}^{d-j}(-1)^{i}{d-i \choose j}f_i, & 0\leq j \leq
\lfloor\frac{d}{2}\rfloor, \nonumber \\
\sum_{i=-1}^{j-1}(-1)^{d-i}{d-i+1 \choose d-j+2}\tilde{f}_i  =
\sum_{i=-1}^{d-j+1}(-1)^{i}{d-i+1 \choose j}\tilde{f}_i, & 0\leq j
\leq \lfloor\frac{d+1}{2}\rfloor.
\end{eqnarray}

Substituting (1) in (2), we get a system of $\lfloor\frac{d}{2}
\rfloor + \lfloor\frac{d+1}{2} \rfloor + 2 = d+2$ independent
linear equations in the $(d-k+1) + (k+1) = d+2$ unknowns $f_k,
\dots, f_d$, $\tilde{f}_{d-k+1}, \dots, \tilde{f}_{d+1}$. Solving
these equations, we get the result (in terms of $f_0, \dots,
f_{k-1}$, which are regarded as ``known" quantities in this
calculation). Notice that this calculation shows that
$\tilde{f}_j$ is given by the same formula as $f_j$ (with $d+1$
in place of $d$ and $\tilde{f}_i = f_i + f_{i-1}$ in place of
$f_i$). This is no surprise: putting $\widetilde{B} = x\ast B$,
one sees that $\widetilde{B}$ is a $(d+2)$-ball with $\partial
\widetilde{B} = \widetilde{X}$ and ${\rm skel}_{d + 1 -
k}(\widetilde{B}) = {\rm skel}_{d+1-k}(\widetilde{X})$. Thus,
$\widetilde{X}$ is also a $k$-stacked sphere. \hfill $\Box$

\bigskip

Now we are ready to state the generalized lower bound
conjecture\,:

\begin{conj} {\bf (GLBC).} \label{GLBC}
For $d \geq 2k+1$, the face-vector $(f_0, \dots, f_d)$ of any
triangulated $d$-sphere $X$ satisfies
$$
f_j \geq \left\{\begin{array}{l} {\displaystyle \sum_{i = -1}^{k
- 1} (-1)^{k-i+1}{j-i-1 \choose j-k}{d -i + 1  \choose j-i}f_{i}},
\hspace{1cm} \mbox{if } ~~ k \leq j \leq d-k, \\
{\displaystyle \sum_{i = -1}^{k - 1} (-1)^{k - i + 1} \left[
{j-i-1 \choose j-k}{d -i + 1  \choose j-i} -{k \choose d-j+1}
{d-i \choose d-k+1} \right.} \\
\mbox{} \hspace{27mm} + {\displaystyle \sum_{l = d-j}^{k - 1}
(-1)^{k - l} \left. {l \choose d-j}{d -i \choose
d-l+1}\right]f_{i}}, \hspace{3mm} \mbox{if } ~~ d-k+1 \leq j \leq
d.
\end{array}
\right.
$$
Equality holds here for some $j$ if and only if $X$ is a
$k$-stacked $d$-sphere.
\end{conj}

\begin{remark}$\!\!\!${\bf .}
{\rm The $k =1$ case of this conjecture is precisely the LBT (for
triangulated spheres). The $j=k$ case of this conjecture was
first stated by McMullen and Walkup \cite{mw} for the smaller
class of polytopal spheres (i.e., boundary complexes of convex
$(d+1)$-polytopes). Note that, when $X$ is a combinatorial
sphere, all its vertex links are spheres, so that using the $j
=k$ case of the conjecture (if settled), one may deduce the
general case by an obvious extension of the M-P-W reduction.

However, note that the vertex links of triangulated spheres need
not be simply connected. (Bj\"{o}rner and Lutz \cite{bl} have
constructed a 16-vertex triangulation $\Sigma^{\,3}_{16}$ of the
Poincar\'{e} homology 3-sphere. Then $S^{\,1}_3 \ast
\Sigma^{\,3}_{16}$ is an example of a triangulated 5-sphere some
of whose vertex-links are not simply connected. Note that the
face-vector of $\Sigma^{\,3}_{16}$ is $(16, 106, 180, 90)$, and
hence the face-vector of the triangulated $5$-sphere $S^{\,1}_3
\ast \Sigma^{\,3}_{16}$ is $(19, 157, 546, 948, 810, 270)$, which
does satisfy Conjecture \ref{GLBC} with $d=5$, $k=2$.) Moreover,
the cases of larger $j$ (the case $j=d$, for instance) of the
conjecture may be easier to settle. In \cite{st1}, Stanley proved
the inequality in Conjecture \ref{GLBC} for polytopal spheres (in
the case $j=k$, but as the vertex links of polytopal spheres are
again polytopal, this settles the inequalities for all $j$).
However, even for polytopal spheres, the case of equality remains
unsolved. It has been suggested that Conjecture \ref{GLBC} holds
for all simply connected triangulated manifolds.}
\end{remark}

We end with a conjecture on non-simply connected triangulated
manifolds.

\begin{conj} {\bf (LBC for the non-simply connected manifolds).}
\label{LBC-NSC} For $d\geq 3$, the face-vector of any connected
and non-simply connected triangulated $d$-manifold $X$ satisfies

$$
f_j(X) \geq \left\{\begin{array}{ll} {d+1 \choose j}f_{0}(X),
& \mbox{if } ~ 1 \leq j < d, \\[1mm]
df_{0}(X), & \mbox{if } ~ j=d.
\end{array}
\right.
$$
Equality holds here for some $j$ if and only if $X$ is obtained
from a stacked $d$-sphere by an elementary handle addition.
\end{conj}

\begin{remark}$\!\!\!${\bf .}
{\rm Notice that Conjecture $\ref{LBC-NSC}$ would imply, in
particular, that the face-vector of any connected and non-simply
connected manifold of dimension $d\geq 3$ must satisfy ${f_0
\choose 2} \geq f_1 \geq (d+1)f_0$, so that any such
triangulation requires $f_0\geq 2d+3$ vertices, and the
triangulation must be 2-neighbourly when $f_0 = 2d+3$. Indeed, in
\cite{bd8}, we proved that any non-simply connected triangulated
$d$-manifold requires at least $2d+3$ vertices, and there is a
unique such $(2d+3)$-vertex triangulated $d$-manifold for $d\geq
3$. It is 2-neighbourly, and does arise from a stacked sphere by
an elementary handle addition. Thus, the main theorem of
\cite{bd8} would be a simple consequence of Conjecture
\ref{LBC-NSC}. The special case $f_0 = 2d+4$ of this conjecture
was posed in \cite{bd8}. In \cite{wa}, Walkup proved that this
conjecture holds for $d = 3$. }
\end{remark}

\bigskip

\noindent {\bf Acknowledgement\,:} The authors are thankful to
Siddhartha Gadgil and Vishwambhar Pati for useful conversations.

\end{document}